\documentclass[a4paper,11pt]{article}
\usepackage{color}

\def\Prg{{\hbox{\rm Prog}}}
\def\Gap{{\hbox{\rm Gap}}}

\def\mones{{\mathbf{1}}}

\def\Opt{\hbox{\rm Opt}}

\def\Prox{{\hbox{}\rm Prox}}
\def\bS{{\mathbf{S}}}
\def\Tr{\trace}
\def\cI{{\cal I}}
\def\cP{{\cal P}}

\def\bR{{\mathbf{R}}}

\def\bE{{\mathbf{E}}}

\def\three{3}

\def\sign{{\hbox{\rm sign}}}

\def\cL{{\cal L}}

\def\three2{2}
\def\cA{{\cal A}}

\def\Argmax{\mathop{\hbox{\rm Argmax}}}
\def\argmin{\mathop{\hbox{\rm argmin}}}

\newtheorem{lemma}{Lemma}[section]
\newtheorem{corollary}{Corollary}[section]

\newtheorem{proposition}{Proposition}[section]

\newtheorem{remark}{Remark}[section]

\newcommand{\qed}{\hfill\hbox{\vrule\vbox{\hrule\phantom{o}\hrule}\vrule}\par\noindent \medskip}

\def\cF{{\cal F}}
\def\cl{{\hbox{\rm cl}\,}}
\def\chapI{Chapter 5}
\def\chapII{Chapter 6}
\def\Tr{\hbox{\rm Tr}}
\def\cg{\hbox{\scriptsize\rm SCG}}
\def\md{\hbox{\scriptsize\rm MD}}
\newcommand{\ese}{\end{eqnarray*}}
\newcommand{\bse}{\begin{eqnarray*}}
\newcommand{\be}{\begin{eqnarray}}
\newcommand{\ee}[1]{\label{#1}\end{eqnarray}}
\newcommand{\nn}{\nonumber \\}
\newcommand{\rf}[1]{~(\ref{#1})}

\begin{document}
\title{Dual subgradient algorithms for large-scale nonsmooth learning problems}
\author{
Bruce Cox\thanks{US Air Force} \and Anatoli Juditsky\thanks{LJK,
Universit\'e J. Fourier, B.P. 53, 38041 Grenoble
Cedex 9, France, {\tt Anatoli.Juditsky@imag.fr}}
\and Arkadi Nemirovski\thanks{Georgia Institute
 of Technology, Atlanta, Georgia
30332, USA, {\tt nemirovs@isye.gatech.edu}\newline
Research of
the third author was supported by the ONR grant N000140811104 and the NSF grants DMS 0914785, CMMI 1232623.}
}
\maketitle
\begin{abstract}
``Classical'' First Order (FO) algorithms of convex optimization, such as Mirror Descent algorithm or Nesterov's optimal algorithm of smooth convex optimization, are well known to have optimal (theoretical) complexity estimates which do not depend on the problem dimension. However, to attain the optimality, the domain of the problem should admit a ``good proximal setup''. The latter essentially means that 1) the problem domain should satisfy certain geometric conditions of ``favorable geometry'', and 2) the practical use of these methods is conditioned by our ability to compute at a moderate cost {\em proximal transformation} at each iteration. More often than not these two conditions are satisfied in optimization problems arising in computational learning, what explains why proximal type FO methods recently became methods of choice when solving various learning problems. Yet, they  meet their limits in several important problems such as multi-task learning with large number of tasks, where the problem domain does not exhibit favorable geometry, and learning and matrix completion problems with nuclear norm constraint, when the numerical cost of computing proximal transformation becomes prohibitive in large-scale problems.
\par
We propose a novel approach to solving nonsmooth optimization problems arising in learning applications where Fenchel-type representation of the objective function is available. The approach is based on applying FO algorithms to the dual problem and using the {\em accuracy certificates} supplied by the method to recover the primal solution. While suboptimal in terms of accuracy guaranties, the proposed approach does not rely upon ``good proximal setup'' for the primal problem but requires the problem domain to admit a {\em Linear Optimization oracle} -- the ability to efficiently maximize a linear form on the domain of the primal problem.
\end{abstract}
\section{Introduction}\label{sect:intro}
\paragraph{Motivation and background.} The problem of interest in this paper is a convex optimization problem in the form
$$
\Opt(P)=\max_{x\in X} f_*(x)\eqno{(P)}
$$
where $X$ is a nonempty closed and bounded subset of Euclidean space $E_x$, and $f_*$ is concave and Lipschitz continuous function on $X$. We are interested in the situation where the sizes of the problem put it beyond the ``practical grasp'' of polynomial time interior point methods with their rather computationally expensive in the large scale case iterations. In this case the methods of choice are the First Order (FO) optimization techniques. The state of the art of these techniques can be briefly summarized as follows:
\\
$\bullet$ The most standard FO approach to ($P$) requires to provide $X,\;E_x$ with {\sl proximal setup}  $\|\cdot\|,\omega(\cdot)$, that is to equip the space $E_x$ with a norm $\|\cdot\|$, and the domain $X$ of the problem -- with a {\sl strongly convex, modulus 1, w.r.t. $\|\cdot\|$ distance-generating function} (d.-g.f.) $\omega(\cdot)$ with variation over $X$ bounded by some $\Omega^2_X$.
    After such a setup is fixed, generating an $\epsilon$-solution to the problem (i.e., a point $x_\epsilon\in X$ satisfying $\Opt(P)-f_*(x_\epsilon)\leq\epsilon$) costs at most $N(\epsilon)$ {\sl steps}, where
    \begin{itemize}
    \item $N(\epsilon)=O(1){\Omega_X^2 L^2\over \epsilon^2}$ in the nonsmooth case, where $f_*$ is Lipschitz continuous, with constant $L$ w.r.t. $\|\cdot\|$ (Mirror Descent (MD) algorithm, see, e.g., \cite[\chapI]{MLO}), and
    \item $N(\epsilon)=O(1){\Omega_X\sqrt{\cL\over\epsilon}}$ in the smooth case, where $f_*$ possesses Lipschitz continuous, with constant $\cL$, gradient: $\|f_*^\prime(x)-f_*^\prime(x')\|_*\leq \cL\|x-x'\|$, where $\|\cdot\|_*$ is the norm conjugate to $\|\cdot\|$ (Nesterov's algorithm for smooth convex optimization, see, e.g., \cite{nesterov1983,NesSmooth}).
Here and in the sequel $O(1)$'s stand for positive absolute constants. Note that in the large scale case, these convergence rates are the best allowed under circumstances by Information-Based Complexity Theory (for details, see,  e.g., \cite{NYu,nesterov_book}).
\end{itemize}
    A {\sl step} of a FO method essentially reduces to computing $f_*$, $f_*^\prime$ at a point and computing { {\sl prox-mapping} (``Bregman projection'') $(x\in X,\xi\in\ E_x)\mapsto \argmin_{u\in X}\left\{\omega(u)+\langle\xi-\omega'(x).u\rangle\right\}$, construction originating from J.-J. Moreau \cite{Mor62,Mor65} \label{Moreau1} and L. Bregman \cite{Bregman67}.}
    \\
$\bullet$ A different way of processing  ($P$) by FO methods, originating in the breakthrough paper of Nesterov \cite{NesSmooth}, is to use {\sl Fenchel-type representation of $f_*$}:
\begin{equation}\label{Fenchel}
f_*(x)=\min\limits_{y\in Y} \left[F(x,y):=\langle x,Ay+a\rangle +\psi(y)\right],
\end{equation}
where $Y$ is a closed and bounded subset of Euclidean space $E_y$, $A:\;E_y\mapsto E_x$ is a linear mapping, and $\psi(y)$ is a convex function.  Representations of this type are readily available for a wide family of ``well-structured'' nonsmooth objectives $f_*$; moreover, usually we can make $\psi$ to possess Lipschitz continuous gradient or even to be linear (for instructive examples, see, e.g., \cite{NesSmooth,Nesterov-dual} or \cite[\chapII]{MLO}). Whenever this is the case, and given proximal setups $(\|\cdot\|_x,\omega_x(\cdot))$ and $(\|\cdot\|_y,\omega_y(\cdot))$ for $(E_x,X)$ and for $(E_y,Y)$, proximal type algorithms like Nesterov's smoothing \cite{NesSmooth} and Dual Extrapolation \cite{Nesterov-dual}, or Mirror Prox \cite{nemirovski2004prox} allow to  find  an $\epsilon$-solution to ($P$) in  $O(1/\epsilon)$ proximal type iterations, the factors hidden in $O(\cdot)$ being explicitly given functions of the variations $\Omega_X^2$ of  $\omega_x$ and $\Omega_Y^2$ of $\omega_y$
and of the partial Lipschitz constants of $\nabla F$ w.r.t. $x$ and to $y$.
\par
Clearly, to be practical, methods of the outlined type should rely on ``good'' proximal setups -- those resulting in ``moderate'' values of $\Omega_{X}$ {\em and} $\Omega_{Y}$ and not too difficult to compute prox-mappings, associated with $\omega_X$ {\em and} $\omega_Y$. This is indeed the case for domains $X$  arising in numerous applications (for instructive examples, see, e.g., \cite[\chapI]{MLO}). The question addressed in this paper is what to do when one of the domains, namely, $X$ does not admit a ``good'' proximal setup. Here are two instructive examples:
\paragraph{A.} $X$ is the unit ball of the nuclear norm $\|\sigma(\cdot)\|_1$ in the space $\bR^{p\times q}$ of $p\times q$ matrices (from now on, for a $p\times q$ matrix $x$, $\sigma(x)=[\sigma_1(x);...;\sigma_{\min[p,q]}(x)]$ denotes the vector comprised by singular values of $x$ taken in the non-ascending order). This domain arises in various low-rank-oriented problems of matrix recovery. In this case, $X$ does admit a proximal setup with $\Omega_{X}=O(1)\sqrt{\ln(pq)}$. However, computing prox-mapping involves full singular value decomposition (SVD) of a $p\times q$ matrix and becomes prohibitively time-consuming when $p,q$ are in the range of tens of thousand. Note that this hardly is a shortcoming of the existing proximal setups, since already computing the nuclear norm (that is, checking the inclusion $x\in X$) requires computing the SVD of $x$.
\paragraph{B.} $X$ is a high-dimensional box -- the unit ball of the $\|\cdot\|_\infty$-norm of $\bR^m$ with large $m$, or, more generally, the unit ball of the $\ell_\infty/\ell_2$ norm $\|x\|_{\infty|2}=\max_{1\leq i\leq m}\|x^i\|_2$, where $x=[x^1;...;x^m]\in E=\bR^{n_1}\times...\times\bR^{n_m}$. Here it is easy to point out a proximal setup with an easy-to-compute prox mapping (e.g., the {\sl Euclidean setup} $\|\cdot\|=\|\cdot\|_2$, $\omega(x)={1\over 2}\langle x,x\rangle$). However, it is easily seen that whenever $\|\cdot\|$ is a norm satisfying $\|e_j\|\geq1$ for all basic orths $e_j$, \footnote{This is a natural normalization: indeed, $\|e_j\|\ll1$ means that $j$-th coordinate of $x$ has a large Lipschitz constant w.r.t. $\|\cdot\|$, in spite of the fact that this coordinate is ``perfectly well behaved'' on $X$ -- its variation on the set is just 2.} one has $\Omega_{X}\geq O(1)\sqrt{m}$, that is, the (theoretical) performance of  ``natural'' proximal FO methods deteriorates rapidly as $m$ grows.
\medskip
\par
Note that whenever a prox-mapping associated with $X$ is ``easy to compute,'' it is equally easy to maximize over $X$ a linear form (since $\Prox_x(-t\xi)$ converges to the maximizer of $\langle\xi,x\rangle$ over $X$ as $t\to \infty$). In such a case, we have at our disposal an efficient {\sl Linear Optimization (LO) oracle} -- a routine which, given on input a linear form $\xi$, returns a point $x_X(\xi)\in\Argmax_{x\in X}\langle\xi,x\rangle$. This conclusion, however, cannot be reversed --  our abilities to maximize, at a reasonable cost, linear functionals over $X$ does not imply the possibility to compute a prox-mapping at a comparable cost. For example, when $X\in\bR^{p\times q}$ is the unit ball of the nuclear norm, maximizing a linear function $\langle\xi,x\rangle=\Tr(\xi x^T)$ over $X$ requires finding the largest singular value of a $p\times q$ matrix and associated left singular vector.  For large $p$ and $q$, solving the latter problem is by orders of magnitude cheaper than computing full SVD of a $p\times q$ matrix. This and similar examples motivate the interest, especially in Machine Learning community, in optimization techniques solving ($P$) via an LO oracle for $X$. In particular, the only ``classical'' technique of this type -- the {\sl Conditional Gradient (CG)} algorithm going back to Frank and Wolfe \cite{Frank1956Algorithm} -- has attracted much attention recently. In the  setting of CG method it is assumed that $f$ is smooth (with H\"{o}lder continuous gradient), and the standard result here (which is widely known, see, e.g., \cite{Dem:Rub:1970,Dunn78,Pshe:1994}) is the following.
\begin{proposition}\label{propredgrad} Let $X$ be a closed and bounded convex set in a Euclidean space $E_x$ such that $X$ linearly spans $E_x$. Assume that we are given  an LO oracle for $X$, and let $f_*$ be a concave continuously differentiable function on $X$ such that for some $\cL<\infty$ and $q\in(1,2]$ one has
\begin{equation}\label{suchthat}
\forall x,x'\in X: f_*(x')\geq f_*(x)+\langle f_*'(x),x'-x\rangle -{1\over q}\cL\|x'-x\|_X^q,
\end{equation}
where $\|\cdot\|_X$ is the norm on $E_x$ with the unit ball ${1\over 2}[X-X]$. Consider a  recurrence of the form
\begin{equation}\label{eqreq}
\begin{array}{lrcl}
&x_t&\mapsto&x_{t+1}\in X: f_*(x_{t+1})\geq f_*\left(x_t+{2\over t+1}[x_X(f_*^\prime(x_t))-x_t]\right),\,t=1,2,...,\\
\end{array}
\end{equation}
where $x_X(\xi)\in\Argmax_{x\in X} \langle\xi,x\rangle$ and $x_1\in X$. Then for all $t=2,3,...$ one has
\begin{equation}\label{then}
\epsilon_{t}:=\max_{x\in X} f_*(x)-f_*(x_{t})\leq {2^{2q}\over q(3-q)}\cdot {\cL\over (t+1)^{q-1}}.
\end{equation}
\end{proposition}
\paragraph{Contents of this paper.}  Assuming an LO oracle for $X$ available, the major limitation in solving ($P$) by the Conditional Gradient method is the requirement for problem objective $f_*$ to be smooth (otherwise, there are no particular requirements to the problem geometry). What to do if this requirement is not satisfied? In this paper, we investigate two simple options for processing this case, based on Fenchel-type representation (\ref{Fenchel}) of $f_*$ {
{\sl which we assume to be available}. Fenchel-type representations are utilized in various ways by a wide spectrum of duality-based convex optimization algorithms (see, e.g., \cite{combettes2010} and references therein). Here we primarily focus on ``nonsmooth'' case, where the representation involves a Lipschitz continuous convex function $\psi$ given by a First Order oracle.} Besides this, we assume that  $Y$ (but not $X$!) {\sl does admit a proximal setup} $(\|\cdot\|_y,\omega_y(\cdot))$. In this case, we can pass from the problem of interest ($P$) to its dual
$$
\Opt(D)=\min\limits_{y\in Y}\left[f(y):=\max\limits_{x\in X} F(x,y)\right]. \eqno{(D)}
$$
Clearly, the LO oracle for $X$ along with the FO oracle for $\psi$ provide a FO oracle for $(D)$:
$$
f(y)=\langle x(y),Ay+a\rangle +\psi(y),\;\;f'(y)=A^Tx(y)+\psi'(y),\;\; x(y):=x_X(Ay+a).
$$
Since $Y$ admits a proximal setup, this is enough to allow to get an $\epsilon$-solution to $(D)$ in $N(\epsilon)=O(1){L^2\Omega^2_{Y}\over\epsilon^2}$ steps, $L$ being the Lipschitz constant of $f$ w.r.t. $\|\cdot\|_y$. Whatever slow the resulting rate of convergence could look, we shall see in the mean time that there are important applications where this rate seems to be the best known so far. When implementing the outlined scheme, the only nontrivial question is how to recover a good optimal solution to the problem  $(P)$ of actual interest from a good approximate solution to its dual problem $(D)$. The proposed answer to this question stems from the pretty simple at the first glance machinery of {\sl accuracy certificates} proposed recently in \cite{NOR}, and closely related to the work \cite{Nesterov-primal}. The summary of our approach is as follows. When solving $(D)$ by a FO method, we generate {\sl search points} $y_\tau\in Y$ where the subgradients $f'(y_\tau)$ of $f$ are computed; as a byproduct of the latter computation, we have at our disposal the points $x_\tau=x(y_\tau)$. As a result, after $t$ steps we have at our disposal {\sl execution protocol} $y^t=\{y_\tau,f'(y_\tau)\}_{\tau=1}^t$. An {\sl accuracy certificate} associated with this protocol is, by definition, a collection  $\lambda^t=\{\lambda^t_\tau\}_{\tau=1}^t$ of nonnegative weights $\lambda^t_\tau$ summing up to 1: $\sum_{\tau=1}^t\lambda^t_\tau=1$. The {\sl resolution} of the certificate is, by definition, the quantity
$$
\epsilon(y^t,\lambda^t)=\max\limits_{y\in Y}\sum_{\tau=1}^t\lambda^t_\tau\langle f'(y_\tau),y_\tau-y\rangle.
$$
An immediate observation is (see section \ref{sect:maindevelopments}) that setting $\widehat{y}^t=\sum_{\tau=1}^t\lambda^t_\tau y_\tau$, $\widehat{x}^t=\sum_{\tau=1}^t\lambda^t_\tau x_\tau
$, we get a pair of feasible solutions to $(D)$ and to $(P)$ such that
$$
[f(\widehat{y}^t)-\Opt(D)]+[\Opt(P)-f_*(\widehat{x}^t)]\leq\epsilon(y^t,\lambda^t).
$$
Thus, assuming that the FO method in question produces, in addition to search points, accuracy certificates for the resulting execution protocols and that {\sl the resolution of these certificates goes to 0 as $t\to\infty$ at some rate}, we can use the certificates to build feasible approximate solutions to $(D)$ {\sl and to $(P)$} with nonoptimalities, in terms of the objectives of the respective problems, going to 0,  at the same rate, as $t\to\infty$.
\par
The scope of the outlined approach depends on whether we are able to equip known methods of nonsmooth convex minimization with computationally cheap mechanisms for building ``good'' accuracy certificates. The meaning of ``good'' in this context is exactly that the rate of convergence of the corresponding resolution to 0 is identical to the standard efficiency estimates of the methods (e.g., for MD this would mean that $\epsilon(y^t,\lambda^t)\leq O(1)L\Omega_{Y}t^{-1/2}$). \cite{NOR} provides a positive answer to this question for the most attractive academically polynomial time oracle-oriented algorithms for convex optimization, like the Ellipsoid method. These methods, however, usually are poorly suited for large-scale applications. In this paper, we provide a positive answer to the above question for the three most attractive oracle-oriented FO methods for {\sl large-scale} nonsmooth convex optimization known to us. Specifically, we consider
\\
\indent $\bullet$ MD (where accuracy certificates are easy to obtain, see also \cite{Nesterov-primal}),\\
\indent $\bullet$ Full Memory Mirror Descent Level (MDL) method (a Mirror Descent extension of the Bundle-Level method \cite{LNN}; to the best of our knowledge, this extension was not yet described in the literature), and\\
\indent $\bullet$ Non-Euclidean Restricted  Memory Level method (NERML) originating from \cite{NERML}, which we believe
is the most attractive tool for large-scale nonsmooth oracle-based convex optimization. To the best of our knowledge, equipping NERML with accuracy certificates is a novel development.
\par
We also consider a different approach to non-smooth convex optimization over a domain given by LO oracle, approach mimicking Nesterov's smoothing \cite{NesSmooth}. Specifically, assuming, as above, that $f_*$ is given by Fenchel-type representation (\ref{Fenchel}) with $Y$ admitting a proximal setup, we use this setup, {\sl exactly in the same way as in \cite{NesSmooth}}, to approximate $f_*$ by a smooth function which then is minimized by the CG algorithm. Therefore, the only difference with \cite{NesSmooth} is in replacing Nesterov's optimal algorithm  for smooth convex optimization (which requires a good proximal  point setup for $X$) with although slower, but less demanding (just LO oracle for $X$ is enough) CG method. We shall see in the mean time that, unsurprisingly, the theoretical complexity of the two outlined approaches -- ``nonsmooth'' and ``smoothing'' ones -- are essentially the same.
\par
The main body of the paper is organized as follows. In section \ref{sect:maindevelopments}, we develop the components of the approach related to duality and show how an accuracy certificate with small resolution
yields a pair of good approximate solutions to $(P)$ and $(D)$. In section \ref{sect:certificates}, we show how to equip the MD, MDL and NERML algorithms with accuracy certificates. In section \ref{sect:smoothing}, we investigate the outlined smoothing approach. In section \ref{sect:examples}, we consider examples, primarily of Machine Learning origin, where
we prone the usage of the proposed algorithms. Section \ref{sect:numerics} reports some preliminary numerical results.
Some technical proofs are relegated to the appendix.

\section{Duality and accuracy certificates}\label{sect:maindevelopments}
\subsection{Situation}\label{sect:situation}
Let $E_x$ be a Euclidean space,
$X\subset E_x$ be a nonempty closed and bounded convex set equipped with {\sl LO oracle} -- a procedure which, given on input $\xi\in E_x$, returns a maximizer $x_X(\xi)$ of the linear form $\langle \xi,x\rangle$ over $x\in X$. Let $f_*(x)$ be a concave function given by {\sl Fenchel-type representation}:
\begin{equation}\label{fench}
f_*(x)=\min_{y\in Y}\left[\langle x,Ay+a\rangle +\psi(y)\right],
\end{equation}
where $Y$ is a convex compact subset of a Euclidean space $E_y$ and $\psi$ is a Lipschitz continuous convex function on $Y$ given by a First Order oracle.
\par
In the sequel we set
$$
f(y)=\max\limits_{x\in X} \left[\langle x,Ay+a\rangle +\psi(y)\right],
$$
and consider two optimization problems
$$
\begin{array}{rcll}
\Opt(P)&=&\max\limits_{x\in X} f_*(x)&(P)\\
\Opt(D)&=&\min\limits_{y\in Y} f(y)&(D)\\
\end{array}
$$
By the standard saddle point argument, we have $\Opt(P)=\Opt(D)$.
\subsection{Main observation}
Observe that the First Order oracle for $\psi$ along with the LO oracle for $X$ provide a First Order oracle for $(D)$; specifically, the vector field
$$
f^\prime(y)=A^Tx_X(Ay+a)+\psi^\prime(y):\;Y\to E_y,
$$
where $\psi^\prime(y)\in \partial \psi(y)$ is a subgradient field of $f$.
\par
Consider a collection $y^t=\{y_\tau\in Y,f^\prime(y_\tau)\}_{\tau=1}^t$ along with a collection $\lambda^t=\{\lambda_\tau\geq0\}_{\tau=1}^t$ such that $\sum_{\tau=1}^t\lambda_\tau=1$, and let us set
$$
\begin{array}{rcl}
y(y^t,\lambda^t)&=&\sum_{\tau=1}^t\lambda_\tau y_\tau,\\
x(y^t,\lambda^t)&=&\sum_{\tau=1}^t\lambda_\tau x_X(Ay_\tau+a),\\
\epsilon(y^t,\lambda^t)&=&\max\limits_{y\in Y} \sum_{\tau=1}^t\lambda_\tau\langle f^\prime(y_\tau),y_\tau-y\rangle.\\
\end{array}
$$
In the sequel, the components $y_\tau$ of $y^t$ will be the search points generated by a First Order minimization method as applied to $(D)$ at the steps $1,...,t$. We call  $y^t$ the associated {\sl execution protocol}, call a collection $\lambda^t$ of $t$ nonnegative weights summing up to 1  an {\sl accuracy certificate} for this protocol, and refer to the quantity $\epsilon(y^t,\lambda^t)$ as to the {\sl resolution of the certificate $\lambda^t$ at the protocol $y^t$}.\par
Our main observation (cf. \cite{NOR}) is as follows:
\begin{proposition}\label{prop1} Let $y^t$, $\lambda^t$ be as above. Then $\widehat{x}:=x(y^t,\lambda^t)$, $\widehat{y}:=y(y^t,\lambda^t)$ are feasible solutions to problems $(P)$, $(D)$, respectively, and
\begin{equation}\label{eqmain}
f(\widehat{y})-f_*(\widehat{x})=
\left[f(\widehat{y})-\Opt(D)\right]+\left[\Opt(P)-f_*(\widehat{x})\right]
\leq\epsilon(y^t,\lambda^t).
\end{equation}
\end{proposition}
{\bf Proof.} Let  $F(x,y)=\langle x,Ay+a\rangle+\psi(y)$ and $x(y)=x_X(Ay+a)$, so that $f(y)=F(x(y),y)$. Observe that  $f^\prime(y)=F^\prime_y(x(y),y)$, where $F^\prime_y(x,y)$ is a selection of the subdifferential of $F$ w.r.t. $y$, that is, $F^\prime_y(x,y) \in\partial_yF(x,y)$ for all $x\in X$, $y\in Y$.  Setting $x_\tau=x(y_\tau)$, we have for all $y\in Y$:
\begin{eqnarray}
\epsilon(y^t,\lambda^t)&\geq& \sum_{\tau=1}^t\lambda_\tau\langle f^\prime(y_\tau),y_\tau-y\rangle=
\sum_{\tau=1}^t\lambda_\tau\langle F^\prime_y(x_\tau,y_\tau),y_\tau-y\rangle\nn
&\geq& \sum_{\tau=1}^t\lambda_\tau\left[F(x_\tau,y_\tau)-F(x_\tau,y)\right]\hbox{\ [by convexity of $F$ in $y$]}\nn
&=&\sum_{\tau=1}^t\lambda_\tau\left[f(y_\tau)-F(x_\tau,y)\right]\hbox{\ [since $x_\tau=x(y_\tau)$, so that $F(x_\tau,y_\tau)=f(y_\tau)$]}\label{upper3}\\
&\geq&f(\widehat{y})-F(\widehat{x},y)\hbox{\ [by convexity of $f$ and concavity of $F(x,y)$ in $x$]}.\nonumber
\end{eqnarray}
We conclude that \[\epsilon(y^t,\lambda^t)\geq \max_{y\in y}\left[f(\widehat{y})-F(\widehat{x},y)\right]=f(\widehat{y}) -f_*(\widehat{x}).\]
The inclusions $\widehat{x}\in X$, $\widehat{y}\in Y$ are evident. \qed
{\small
\begin{remark}\label{rem1} {\rm In the proof of Proposition \ref{prop1}, the linearity of $F$ w.r.t. $x$ was never used, so that in fact we have proved a more general statement:}\\
Given a concave in $x\in X$ and convex in $y\in Y$ Lipschitz continuous function $F(x,y)$, let us associate with it a convex function $f(y)=\max_{x\in X} F(x,y)$, a concave function $f_*(x)=\min_{y\in Y}F(x,y)$ and problems $(P)$ and $(D)$. Let $F^\prime_y(x,y)$ be a vector field with $F^\prime_y(x,y)\in\partial_yF(x,y)$, so that with $x(y)\in\Argmax_{x\in X} F(x,y)$, the vector $f^\prime(y)=F^\prime_y(x(y),y)$ is a subgradient of $f$ at $y$. Assume that problem $(D)$ associated with $F$ is solved by a FO method using $f^\prime(y)=F^\prime_y(x(y),y)$ which produced execution protocol $y^t$ and accuracy certificate $\lambda^t$. Then setting
\[\widehat{x}=\sum_\tau\lambda_\tau x(y_\tau),\;\; \mbox{and}\;\; \widehat{y}=\sum_\tau\lambda_\tau y_\tau,
\] we ensure {\rm (\ref{eqmain})}.\par
Moreover, let $\delta\geq0$, and let $x_\delta(y)$ be a $\delta$-maximizer of $F(x,y)$ in $x\in X$: for all $y\in Y$,
$$
F(x_\delta(y),y)\geq {\max}_{x\in X} F(x,y)-\delta.
$$
Suppose that $(D)$ is solved by a FO method using approximate subgradients $\tilde{f}^\prime(y)=F^\prime_y(x_\delta(y),y)$, and producing execution protocol $y^t=\{y_\tau,\tilde{f}^\prime(y_\tau)\}_{\tau=1}^t$ and accuracy certificate $\lambda^t$. Then setting $\widehat{x}=\sum_\tau\lambda_\tau x_\delta(y_\tau)$ and $\widehat{y}=\sum_\tau\lambda_\tau y_\tau$, we ensure the $\delta$-relaxed version of {\rm (\ref{eqmain})} -- the relation
$$
f(\widehat{y})-f_*(\widehat{x}) \leq\epsilon(y^t,\lambda^t) + \delta,\,\,\epsilon(y^t,\lambda^t)=\max_{y\in Y} \sum_{\tau=1}^t\lambda^t_\tau\langle \tilde{f}^\prime(y_\tau),y_\tau-y\rangle.
$$
{\rm All we need to extract the ``Moreover'' part of this statement from the proof of Proposition \ref{prop1} is to set $x_\tau=x_\delta(y_\tau)$, to replace $f'(y_\tau)$ with $\tilde{f}^\prime(y_\tau)$   and to replace the equality in \rf{upper3} with the inequality $$\sum_{\tau=1}^t\lambda_\tau\left[F(x_\tau,y_\tau)-F(x_\tau,y)\right]\geq
\sum_{\tau=1}^t\lambda_\tau\left[f(y_\tau)-\delta-F(x_\tau,y)\right].$$}
\end{remark}}
\paragraph{Discussion.} Proposition \ref{prop1} says that whenever we can equip the subsequent execution protocols generated by a FO method,
 as applied to the dual problem $(D)$,  with accuracy certificates,
 we can generate solutions to the primal problem $(P)$ of inaccuracy going to 0 at the same rate as the certificate resolution. In the sequel, we shall point out some ``good'' accuracy certificates for several most attractive FO algorithms for nonsmooth convex minimization.

\section{Accuracy certificates in oracle-oriented methods for large-scale nonsmooth convex optimization}
\label{sect:certificates}
\subsection{Convex minimization with certificates, I: Mirror Descent}
\subsubsection{Proximal setup.}\label{ProximalSetup} As it was mentioned in the introduction, the Mirror Descent (MD) algorithm solving $(D)$ is given by a norm $\|\cdot\|$ on $E_y$ and a distance-generating function (d.-g.f.) $\omega(y):Y\to\bR$ which should be continuous and convex on $Y$, should admit a continuous in $y\in Y^o=\{y\in Y:\partial\omega(y)\neq\emptyset\}$ selection of subdifferentials $\omega'(y)$, and should be strongly convex, modulus 1, w.r.t. $\|\cdot\|$, that is,
$$
\forall y,y'\in Y^o: \langle\omega'(y)-\omega'(y'),y-y'\rangle \geq \|y-y'\|^2.
$$
A proximal setup $(\|\cdot\|,\omega(\cdot))$ for $Y,E_y$ gives rise to several entities, namely,
\begin{itemize}
\item Bregman distance $V_y(z)=\omega(z)-\omega(y)-\langle\omega'(y),z-y\rangle$ ($y\in Y^o, z\in Y$). Due to strong convexity of $\omega$, we have
\begin{equation}\label{wehave22}
\forall (z\in Y, y\in Y^o): V_y(z)\geq {1\over 2}\|z-y\|^2;
\end{equation}
\item $\omega$-center $y_\omega=\argmin_{y\in Y}\omega(y)$ of $Y$ and $\omega$-diameter
\[\Omega=\Omega[Y,\omega(\cdot)]:=\sqrt{2\left[\max_{y\in Y}\omega(y)-\min_{y\in Y}\omega(y)\right]}.
\] Observe that
\begin{equation}\label{AppLem1}
\langle \omega'(y_\omega),y-y_\omega\rangle \geq0,
\end{equation} (see Lemma \ref{AppLemma}), so that
    \begin{equation}\label{Visbounded}
    V_{y_\omega}(z)\leq \omega(z)-\omega(y_\omega)\leq {1\over 2}\Omega^2,\,\forall z\in Y,
     \end{equation}
     which combines with the inequality $V_y(z)\geq {1\over 2}\|z-y\|^2$ to yield the relation
    \begin{equation}\label{diameter}
    \forall y\in Y: \|y-y_\omega\|\leq \Omega;
    \end{equation}
\item prox-mapping
$$
\Prox_y(\xi)=\argmin_{z\in Y}\left[\langle\xi,z\rangle + V_y(z)\right],
$$
where $\xi\in E_y$ and $y\in Y^o$. This mapping takes its values in $Y^o$ and satisfies the relation \cite{Teboulle93}
\begin{equation}\label{identity}
\forall (y\in Y^o,\xi\in E_y,y_+=\Prox_y(\xi)): \langle \xi,y_+-z\rangle \leq V_y(z)-V_{y_+}(z)-V_y(y_+)\,\,\forall z\in Y.
\end{equation}
\end{itemize}
\subsubsection{Mirror Descent algorithm} MD algorithm works with a vector field
 \begin{equation}\label{field}
y\mapsto g(y): Y\to E_y,
\end{equation}
which is {\sl oracle represented}, meaning that we have access to an oracle which, given on input $y\in Y$, returns $g(y)$. From now on we assume that this field is bounded:
\begin{equation}\label{normalization}
\|g(y)\|_*\leq L[g]<\infty,\,\forall y\in Y,
\end{equation}
where $\|\cdot\|_*$ is the norm conjugate to $\|\cdot\|$. The algorithm is the recurrence
$$
y_1=y_\omega;y_\tau\mapsto g_\tau:=g(y_\tau)\mapsto y_{\tau+1}:=\Prox_{y_\tau}(\gamma_\tau g_\tau),\eqno{\hbox{(MD)}}
$$
where $\gamma_\tau>0$ are stepsizes. Let us equip this recurrence with accuracy certificates, setting
\begin{equation}\label{lambdat}
\lambda^t=\left({\sum}_{\tau=1}^t\gamma_\tau\right)^{-1}[\gamma_1;...;\gamma_t].
\end{equation}
\begin{proposition} \label{PropMD} For every $t$, the resolution
$$
\epsilon(y^t,\lambda^t):=\max\limits_{y\in Y}\sum_{\tau=1}^t\lambda_\tau\langle g(y_\tau),y_\tau-y\rangle
$$
of  $\lambda^t$ on the execution protocol $y^t=\{y_\tau,g(y_\tau)\}_{\tau=1}^t$ satisfies the standard MD efficiency estimate
\begin{equation}\label{MDcase}
\epsilon(y^t,\lambda^t) \leq {\Omega^2+\sum_{\tau=1}^t\gamma_\tau^2\|g(y_\tau)\|_*^2\over 2\sum_{\tau=1}^t\gamma_\tau}.
\end{equation}
In particular, if $\gamma_\tau={\gamma(t)\over \|g(y_\tau)\|_*}$, $\gamma(t):={\Omega\over\sqrt{t}}$ for $1\leq\tau\leq t$  \footnote{We assume here that $g_\tau\neq0$ for all $\tau\leq t$. In the opposite case, the situation is trivial: when $g(y_{\tau_*})=0$, for some $\tau*\leq t$, setting $\lambda^t_\tau=0$ for $\tau\neq\tau_*$ and $\lambda^t_{\tau_*}=1$, we ensure that $\epsilon(y^t,\lambda^t)=0$.},
\begin{equation}\label{MDini}
\epsilon(y^t,\lambda^t)\leq {\Omega L[g]\over\sqrt{t}}.
\end{equation}
\end{proposition}
The proof of the proposition follows the lines of the ``classical'' proof in the case when $g(\cdot)$ is the (sub)gradient field of the objective of a convex minimization problem (see, e.g., Proposition 5.1 of \cite{MLO}), and is omitted.
\paragraph{Solving $(P)$ and $(D)$ using MD} In order to solve $(D)$, we apply MD to the vector field $g=f'$.
Assuming that
\begin{equation}\label{Lf}
L_f=\sup_{y\in Y}\{\|f'(y)\|_*\equiv \|A^*x(y)+\psi'(y)\|_*\}<\infty \quad [x(y)=x_X(Ay+a)],
\end{equation}
we can set $L[g]=L_f$. For this setup, Proposition \ref{prop1} implies that the MD accuracy certificate $\lambda^t$, as defined in (\ref{lambdat}), taken together with the MD execution protocol $y^t=\{y_\tau,g(y_\tau)=f'(y_\tau):=A^*\underbrace{x(y_\tau)}_{x_\tau}+\psi'(y_\tau)\}_{\tau=1}^t$, yield
the primal-dual feasible approximate solutions
\begin{equation}\label{eq229}
\widehat{x}^t=\sum_{\tau=1}^t\lambda^t_\tau x_\tau,\,\,\widehat{y}^t=\sum_{\tau=1}^t\lambda^t_\tau y_\tau
\end{equation}
to $(P)$ and to $(D)$
such that
$$
f(\widehat{y}^t)-f_*(\widehat{x}^t)\le \epsilon(y^t,\lambda^t).
$$
Combining this conclusion with Proposition \ref{PropMD}, we arrive at the following result:
\begin{corollary}\label{cor1} In the case of {\rm (\ref{Lf})}, for every $t=1,2,...$ the $t$-step MD with the stepsize policy\footnote{We assume that $f'(y_\tau)\neq0$, for $\tau\leq t$; otherwise, as we remember, the situation is trivial.}
$$
\gamma_\tau={\Omega\over\sqrt{t}\|f'(x_\tau)\|_*},\,1\leq \tau\leq t
$$
as applied to $(D)$ yields feasible approximate solutions $\widehat{x}^t$, $\widehat{y}^t$ to $(P)$, $(D)$ such that
\begin{equation}\label{MDbound}
\begin{array}{c}
[f(\widehat{y}^t)-\Opt(P)]+[\Opt(D)-f_*(\widehat{x}^t)]\leq {\Omega L_f\over\sqrt{t}}\\
\left[\Omega=\Omega[Y,\omega(\cdot)]\right]\\
\end{array}
\end{equation}
In particular, given $\epsilon>0$, it takes at most
\begin{equation}\label{MDefficiency}
t(\epsilon)=\hbox{\rm Ceil} \left({\Omega^2L_f^2\over\epsilon^2}\right)
\end{equation}
steps of the algorithm to ensure that
\begin{equation}\label{accuracy}
[f(\widehat{y}^t)-\Opt(P)]+[\Opt(D)-f_*(\widehat{x}^t)]\leq\epsilon.
\end{equation}
\end{corollary}

\subsection{Convex minimization with certificates, II: Mirror Descent with full memory}\label{sect:MDL}
Algorithm MDL -- Mirror Descent Level method -- is a non-Euclidean version of (a variant of) the Bundle-Level method \cite{LNN}; to the best of our knowledge, this extension was not presented in the literature. Another novelty in what follows is equipping the method with accuracy certificates.\par
MDL is a version of MD with ``full memory'', meaning that the first order information on the objective being minimized is preserved and utilized at subsequent steps, rather than being ``summarized'' in the current iterate, as it is the case for MD. While the guaranteed accuracy bounds for MDL are similar to those for MD, the typical practical behavior of the algorithm is better than that of the ``memoryless'' MD.
\subsubsection{Preliminaries}
MDL with certificates which we are about to describe is aimed at processing an oracle-represented vector field (\ref{field}) satisfying (\ref{normalization}), with the same assumptions on $Y$ and the same proximal setup as in the case of MD.\par
We associate with $y\in Y$ the affine function
$$
h_y(z)=\langle g(y),y-z\rangle\;\quad[\ge f(y)-f(z)\hbox{\ when $g$ is the subgradient field of convex function $f$}],
$$
and with a finite set $S\subset Y$ the family $\cF_S$ of affine functions on $E_y$ which are convex combinations of the functions $h_y(z)$, $y\in S$.
In the sequel, the words ``we have at our disposal a function $h(\cdot)\in \cF_S$'' mean that we know the functions $h_y(\cdot)$, $y\in S$, and nonnegative weights $\lambda_y$, $y\in S$, summing up to 1, such that $h(z)=\sum_{y\in S}\lambda_yh_y(z)$.
\paragraph{The goal} of the algorithm is, given a tolerance $\epsilon>0$, to find a finite set $S\subseteq Y$ and $h\in \cF_S$ such that
\begin{equation}\label{bruce1}
\max\limits_{y\in Y} h(y)\leq\epsilon.
\end{equation}
Note that our target $S$ and $h$ are of the form  $S=\{y_1,...,y_t\}$, $h(y)=\sum_{\tau=1}^t\lambda_\tau \langle g(y_\tau),y_\tau-y\rangle$ with nonnegative $\lambda_\tau$ summing up to 1. In other words, our target is to build an execution protocol $y^t=\{y_\tau,g(y_\tau)\}_{\tau=1}^t$ {\sl and} an associated accuracy certificate $\lambda^t$ such that $\epsilon(y^t,\lambda^t)\leq\epsilon$.
\subsubsection{Construction}
As applied to (\ref{field}), MDL at a step $t=1,2,...$ generates search point $y_t\in Y$ where the value $g(y_t)$ of $g$ is computed;
it provides us with the affine function $h_t(z)=\langle g(y_t),y_t-z\rangle$. Besides this, the method generates finite sets $I_t\subset\{1,...,t\}$ and
$$
S_t=\{y_\tau:\tau\in I_t\}\subset Y.
$$
Steps of the method are split into subsequent {\sl phases} numbered $s=1,2,...$, and every phase is associated with {\sl optimality gap} $\Delta_s\geq 0$.
\par To initialize the method, we set $y_1=y_\omega$, $I_0=\emptyset$ (whence $S_0=\emptyset$ as well), $\Delta_0=+\infty$. \par
At a step $t$
we act as follows:
\begin{itemize}
\item given $y_t$, we compute $g(y_t)$, thus getting $h_t(\cdot)$, and set
$I_t^+=I_{t-1}\cup\{t\}$, $S_t^+=S_{t-1}\cup\{y_t\}$;
\item we solve the auxiliary problem
\begin{equation}\label{auxlevel1}
\epsilon_t=\max_{y\in Y}\min_{\tau\in I_t^+}h_\tau(y)=\max_{y\in Y}\min\limits_{\lambda=\{\lambda_\tau\}_{\tau\in I_t^+}}\left\{\sum_{\tau\in I_t^+}\lambda_\tau h_\tau(y): \lambda\geq0,\sum_{\tau\in I_t^+}\lambda_\tau=1\right\}
\end{equation}
By the von Neumann lemma, an optimal solution to this (auxiliary) problem is associated with nonnegative and summing up to 1  weights $\lambda_\tau^t$, $\tau\in I_t^+$ such that
$$
\epsilon_t=\max_{y\in Y}\sum_{\tau\in I_t^+}\lambda^t_\tau h_\tau(y),
$$
and we assume that as a result of solving (\ref{auxlevel1}), both $\epsilon_t$ and $\lambda^t_\tau$ become known. We set $\lambda^t_\tau=0$ for all $\tau\leq t$ which are not in $I_t^+$, thus getting an accuracy certificate $\lambda^t=[\lambda^t_1;...;\lambda^t_t]$ for the execution protocol $y^t=\{y_\tau,g(y_\tau)\}_{\tau=1}^t$ along with $h^t(\cdot)=\sum_{\tau=1}^t\lambda^t_\tau h_\tau(\cdot)$. Note that by construction
$$
\epsilon(y^t,\lambda^t)=\max\limits_{y\in Y}h^t(y)=\epsilon_t.
$$
If $\epsilon_t\leq\epsilon$, we terminate -- $h(\cdot)=h^t(\cdot)$ satisfies (\ref{bruce1}). Otherwise we proceed as follows:
\item If (case A) $\epsilon_t\leq \gamma\Delta_{s-1}$, $\gamma\in(0,1)$ being method's control parameter, we say that step $t$ starts phase $s$ (e.g., step $t=1$ starts phase 1),  set
     $$
     \Delta_s=\epsilon_t,\,\,I_t=\{\tau:1\leq\tau\leq t:\lambda^t_\tau>0\}\cup\{1\},\,S_t=\{y_\tau:\tau\in I_t\},\,\widehat{y}_t=y_\omega\\
     $$
    otherwise (case B) we set
    $$
    I_t=I_t^+,\,S_t=S_t^+,\,\widehat{y}_t=y_t.
    $$
Note that in both cases we have
\begin{equation}\label{bothcases}
\epsilon_t=\max_{y\in Y}\min_{\tau\in I_t} h_\tau(y)
\end{equation}
\item Finally, we define {\sl $t$-th level} as $\ell_t=\gamma\epsilon_t$ and associate with this quantity the {\sl level set} $U_t=\{y\in Y: h_\tau(y)\geq \ell_t\;\forall \tau\in I_t\}$, specify $y_{t+1}$ as the $\omega$-projection of $\widehat{y}_t$ on $U_t$:
\begin{equation}\label{auxlevel2}
y_{t+1}=\argmin_{y\in U_t}\left[\omega(y)-\langle\omega'(\widehat{y}_t),y\rangle\right]
\end{equation}
and loop to step $t+1$.
\end{itemize}
\subsubsection{Efficiency estimate}
\begin{proposition}\label{PropLevel}
Given on input a target tolerance $\epsilon>0$, the MDL algorithm terminates after finitely many steps, with the output $y^t=\{Y\ni y_\tau,g(y_\tau)\}_{\tau=1}^t$, $\lambda^t=\{\lambda^t_\tau\geq0\}_{\tau=1}^t$,  $\sum_\tau \lambda^t_\tau=1$ such that
\begin{equation}\label{ressmall}
\epsilon(y^t,\lambda^t)\leq\epsilon.
\end{equation}
The number of steps of the algorithm does not exceed
\begin{equation}\label{numberofstepsLevel}
N={2\Omega^2L^2[g]\over\gamma^4(1-\gamma^2)\epsilon^2}+1 \qquad[\Omega=\Omega[Y,\omega(\cdot)]].
\end{equation}
\end{proposition}
For proof, see Section \ref{Appendix.PropLevel}.
\begin{remark}\label{remLevel} {\rm Assume that $\omega(\cdot)$ is continuously differentiable on the entire $Y$, so that the quantity
$$
\Omega^+=\Omega^+[Y,\omega]:=\max\limits_{y,z\in Y} V_y(z)
$$ is finite. From the proof of Proposition \ref{PropLevel} it follows immediately  that one can substitute the rule ``$\widehat{y}_t=y_\omega$ when $t$ starts a phase and $\widehat{y}_t=y_t$ otherwise'' with a simpler one ``$\widehat{y}_t=y_t$ for all $t$,'' at the price of replacing $\Omega$ in (\ref{numberofstepsLevel}) with $\Omega^+$.}
\end{remark}
\paragraph{Solving $(P)$ and $(D)$ via MDL} is completely similar to the case of MD: given a desired tolerance $\epsilon>0$, one applies MDL to the vector field $g(y)=f'(y)$ until the target (\ref{bruce1}) is satisfied.  Assuming (\ref{Lf}), we can set $L[g]=L_f$, so that by Proposition \ref{PropLevel} our target will be achieved in
\begin{equation}\label{MDLbound}
t(\epsilon)\leq \hbox{\rm Ceil}\left({\Omega^2L_f^2\over\gamma^4(1-\gamma^2)\epsilon^2}+1\right)\quad [\Omega=\Omega[Y,\omega(\cdot)]]
\end{equation}
steps, with $L_f$ given by (\ref{Lf}). Assuming that the target is attained at a step $t$,
we have at our disposal the execution protocol $y^t=\{y_\tau,f'(y_\tau)\}_{\tau=1}^t$ along with the accuracy certificate $\lambda^t=\{\lambda^t_\tau\}$ such that $\epsilon(y^t,\lambda^t)\leq\epsilon$  (by the same Proposition \ref{PropLevel}). Therefore, specifying $\widehat{x}^t$, $\widehat{y}^t$ according to (\ref{eq229}) and invoking Proposition \ref{prop1}, we ensure (\ref{accuracy}).  Note that the complexity $t=t(\epsilon)$ of finding these solutions, as given by (\ref{MDLbound}), is completely similar to the complexity bound (\ref{MDefficiency}) of MD.
\subsection{Convex minimization with certificates, III: restricted memory Mirror Descent}\label{NERML}
The fact that the number of linear functions
$h_\tau(\cdot)$ involved into the auxiliary problems (\ref{auxlevel1}), (\ref{auxlevel2}) (and thus computational complexity of these problems) grows as the algorithm proceeds is a serious shortcoming of MDL from the computational viewpoint.
NERML (Non-Euclidean Restricted Memory Level Method) algorithm, which originates from \cite{NERML} is a version of MD ``with restricted memory''. In this algorithm the number of affine functions $h_\tau(\cdot)$ participating in (\ref{auxlevel1}), (\ref{auxlevel2}) never exceeds $m+1$, where $m$ is a control parameter which can be set to any desired positive integer value. The original NERML algorithm, however, was not equipped with accuracy certificates, and our goal here is to correct this omission.
\subsubsection{Construction}Same as MDL, NERML processes an oracle-represented vector field (\ref{field}) satisfying the boundedness condition (\ref{normalization}), with the ultimate goal to ensure (\ref{bruce1}). The setup for the algorithm is identical to that for MDL.

The  algorithm builds search sequence $y_1\in Y,y_2\in Y,...$ along with the sets $S_\tau=\{y_1,...,y_\tau\}$, according to the following rules:\\
{\bf A. Initialization.} We set $y_1=y_\omega:=\argmin_{y\in Y}\omega(y)$, compute $g(y_1)$  and set
$f_1=\max\limits_{y\in Y}h_{y_1}(y)$. We clearly
have $f_1\geq0$.
\begin{itemize}
\item In the case of $f_1=0$, we terminate and output $h(\cdot)=h_{y_1}(\cdot)\in\cF_{S_1}$, thus ensuring (\ref{bruce1}) with $\epsilon=0$.
\item When $f_1>0$, we proceed. Our subsequent actions are split into {\sl phases} indexed with $s=1,2,...$.
\end{itemize}
{\bf B. Phase $s=1,2,...$}  At the beginning of phase $s$, we have at our disposal
\begin{itemize}
\item the set $S^s=\{y_1,...,y_{t_s}\}\subset Y$ of already built search points, and
\item an affine function $h^s(\cdot)\in \cF_{S^s}$ along with the real $f_s:=\max\limits_{y\in Y} h^s(y)\in(0,f_1]$.
\end{itemize}
We define the {\sl level} $\ell_s$ of phase $s$ as
$$
\ell_s=\gamma f_s,
$$
where $\gamma\in(0,1)$ is a control parameter of the method. Note that $\ell_s>0$ due to $f_s>0$.
\par To save notation, we denote the search points generated at phase $s$ as $u_1,u_2,...$, so that
$y_{t_s+\tau}=u_{\tau}$, $\tau=1,2,...$.
\\
{\bf B.1. Initializing phase $s$.} We somehow choose collection of $m$ functions $h_{0,j}^s(\cdot)\in\cF_{S^s}$, $1\leq j\leq m$, such that the set $$Y^s_0=\cl\{y\in Y:
    h_{0,j}^s(y)>\ell_s,1\leq j\leq m\}$$ is nonempty (here  a positive integer $m$ is a control parameter of the method).\footnote{Note that to ensure the nonemptiness of $Y^s_0$, it suffices to set $h_{0,j}^s(\cdot)=h^s(\cdot)$, so that $h_{0,j}(y)>\ell_s$ for $y\in\Argmax_Yh^s(\cdot)$; recall that $f_s=\max_{y\in Y} h^s(y)>0$.} We set
$$
u_1=y_\omega.
$$
{\bf B.2. Step $\tau=1,2,...$ of phase $s$:}\\
{\bf B.2.1.} At the beginning of step $\tau$, we have at our disposal
\begin{enumerate}
\item the set $S^s_{\tau-1}$ of all previous search points;
\item a collection of functions $\{h^s_{\tau-1,j}(\cdot)\in\cF_{S^s_{\tau-1}}\}_{j=1}^m$  such that the set $$Y^s_{\tau-1}=\cl\{x\in Y:
    h_{\tau-1,j}^s(x)>\ell_s,1\leq j\leq m\}$$ is nonempty,
\item current search point $u_\tau\in Y^s_{\tau-1}$ such that
$$
u_\tau=\argmin_{y\in Y^s_{\tau-1}} \omega(y)\eqno{(\Pi^s_\tau)}
$$
Note that this relation is trivially true when $\tau=1$.
\end{enumerate}
{\bf B.2.2.} Our actions at step $\tau$ are as follows.\\
{\bf B.2.2.1.} We compute $g(u_\tau)$ and set
$$
h_{\tau-1,m+1}(y)=\langle g(u_\tau),u_\tau-y\rangle.
$$
{\bf B.2.2.2.} We solve the auxiliary problem
\begin{equation}\label{bruce10}
\Opt=\max\limits_{y\in Y} \min\limits_{1\leq j\leq m+1}h_{\tau-1,j}(y).
\end{equation}
Note that
$$
\begin{array}{rcl}
\Opt&=&\max\limits_{y\in Y}\min\limits_{\lambda_j\geq0,\sum_j\lambda_j=1}\sum\limits_{j=1}^{m+1}
\lambda_jh^s_{\tau-1,j}(y)=\min\limits_{\lambda_j\geq0,\sum_j\lambda_j=1}\max\limits_{y\in Y}
\sum\limits_{j=1}^{m+1}
\lambda_jh^s_{\tau-1,j}(y)\\
&=&\max\limits_{y\in Y} \sum_{j=1}^{m+1}\lambda_j^\tau h^s_{\tau-1,j}(y),\\
\end{array}
$$
$\lambda_j^\tau\geq0$ and $\sum_{j=1}^{m+1}\lambda^\tau_j=1$. We assume that when solving the auxiliary problem, we compute the above weights $\lambda^\tau_j$, and thus have at our disposal the function
$$
h^{s,\tau}(\cdot)=\sum_{j=1}^{m+1}\lambda^\tau_jh^s_{\tau-1,j}(\cdot)\in \cF_{S^s_\tau}
$$
such that
$$
\Opt=\max\limits_{y\in Y}h^{s,\tau}(y).
$$
{\bf B.2.2.3.} {\bf Case A:} If  $\Opt\leq\epsilon$ we terminate and output $h^{s,\tau}(\cdot)\in\cF_{S^s_\tau}$; this function satisfies (\ref{bruce1}).\\
{\bf Case B:}
In case of $\Opt<\ell_s+\theta(f_s-\ell_s)$, where $\theta\in(0,1)$ is method's control parameter, we terminate phase $s$ and start phase $s+1$ by setting $h^{s+1}=h^{s,\tau}$, $f_{s+1}=\Opt$. Note that by construction
$0<f_{s+1}\leq [\gamma+\theta(1-\gamma)] f_s< f_1$, so that we have at our disposal all we need to start phase $s+1$.\\
{\bf Case C:} When neither A nor B takes place, we proceed with phase $s$, specifically, as follows:
\\
{\bf B.2.2.4.} Note that there exists a point $u\in Y$ such that $h^s_{\tau-1,j}(u)\geq \Opt>\ell_s$, so that the set $Y_\tau=\{y\in Y: h^s_{\tau-1,j}(y)\geq\ell_s,\,1\leq j\leq m+1\}$, intersects with the relative interior of $Y$. We specify $u_{\tau+1}$ as
\begin{equation}\label{specify}
    u_{\tau+1}=\argmin_{y\in Y_\tau} \omega(y).
    \end{equation}
    Observe that
    \begin{equation}\label{bruce17}
    u_{\tau+1}\in Y^s_{\tau-1}
    \end{equation}
    due to $Y_\tau\subset Y^s_{\tau-1}$.
\\
{\bf B.2.2.5.} By optimality conditions for (\ref{specify})  (see Lemma \ref{AppLemma}), for certain nonnegative $\mu_j$, $1\leq j\leq m+1$, such that
    $$
    \mu_j[h^s_{\tau-1,j}(u_{\tau+1})-\ell_s]=0,\,1\leq j\leq m+1,
    $$
    the vector
\begin{equation}\label{bruce2}
e:=\omega^\prime(u_{\tau+1})-\sum\limits_{j=1}^{m+1}\mu_j\nabla h^s_{\tau-1,j}(\cdot)
\end{equation}
is such that
\begin{equation}\label{bruce3}
\langle e,y-u_{\tau+1}\rangle \geq 0\,\,\forall y\in Y.
\end{equation}
\par
$\bullet$ In the case of $\mu=\sum_j\mu_j>0$, we set
$$
h^s_{\tau,1}={1\over \mu}\sum\limits_{j=1}^{m+1}\mu_jh^s_{\tau-1,j},
$$
so that
\begin{equation}\label{bruce4}
\begin{array}{llllll}
(a)&h^s_{\tau,1}\in\cF_{S^s_\tau},&
(b)&h^s_{\tau,1}(u_{\tau+1})=\ell_s,&
(c)&\langle \omega^\prime(u_{\tau+1})-\mu\nabla h^s_{\tau,1},y-u_{\tau+1}\rangle \geq0\,\forall y\in Y\\
\end{array}
\end{equation}
We then discard from the collection $\{h^s_{\tau-1,j}(\cdot)\}_{j=1}^{m+1}$ two (arbitrarily chosen) elements and add to $h^s_{\tau,1}$ the remaining $m-1$ elements of the collection, thus getting an $m$-element collection $\{h^s_{\tau,j}\}_{j=1}^m$ of elements of $\cF_{S^s_\tau}$.
\begin{quote}
{\small \begin{remark}\label{remNERML}{\rm  We have ensured that the set $Y^s_\tau=\cl\{y\in Y:h^s_{\tau,j}(y)>\ell_s,1\leq j\leq m\}$ is nonempty (indeed, we clearly have $h^s_{\tau,j}(\widehat{u})>\ell_s$, $1\leq j\leq m$, where $\widehat{u}$ is an optimal solution to (\ref{bruce10})).
Besides this, we also have
($\Pi^s_{\tau+1})$. Indeed, by construction $u_{\tau+1}\in Y_\tau$, meaning that $h^s_{\tau-1,j}(u_{\tau+1})\geq\ell_s$, $1\leq j\leq m+1$. Since $h^s_{\tau,j}$ are convex combinations of the functions  $h^s_{\tau-1,j}$, $1\leq j\leq m+1$, it follows that $u_{\tau+1}\in Y^s_{\tau}$. Further,  (\ref{bruce4}.$b$) and (\ref{bruce4}.$c$) imply that
$
u_{\tau+1}=\argmin_y\left\{\omega(y):\,y\in Y, \,h^s_{\tau,1}(y)\geq\ell_s\right\},
$ and the right hand side set clearly contains $Y^s_{\tau}$. We conclude that $u_{\tau+1}$ indeed is the minimizer of $\omega(\cdot)$ on $Y^s_{\tau}$.
}
\end{remark}}
\end{quote}
$\quad$
$\bullet$ In the case of $\mu=0$, (\ref{bruce2}) -- (\ref{bruce3}) say that $u_{\tau+1}$ is a minimizer of $\omega(\cdot)$ on $Y$. In this case, we discard from the collection $\{h^s_{\tau-1,j}\}_{j=1}^{m+1}$ one (arbitrarily chosen) element, thus getting the $m$-element collection $\{h^s_{\tau,j}\}_{j=1}^m$. Here, by exactly the same argument as above,  the set
$Y^s_{\tau}:=\cl\{y\in Y: h^s_{\tau,j}(y)>\ell_s\}$ is nonempty and contains $u_{\tau+1}$, and, of course, ($\Pi^s_{\tau}$) holds true (since $u_{\tau+1}$ minimizes $\omega(\cdot)$ on the entire $Y$).
\par
In both cases (those of $\mu>0$ and of $\mu=0$), we have built the data required to start step $\tau+1$ of phase $s$, and we proceed to this step.
\\
The description of the algorithm is completed.
\begin{remark}\label{remA}{\rm Same as MDL, the outlined algorithm requires solving at every step two nontrivial auxiliary optimization problems -- (\ref{bruce10}) and (\ref{specify}). It is explained in \cite{NERML} that these problems are relatively easy, provided that $m$ is moderate (note that this parameter is under our full control) and $Y$ and $\omega$ are ``simple and fit each other,'' meaning that we can easily solve problems of the form
$$
\min_{x\in Y} \left[\omega(x)+\langle a,x\rangle\right]\eqno{(*)}
$$
(that is, our proximal setup for $Y$ results in easy-to-compute prox-mapping). }
\end{remark}
\begin{remark}\label{remB}{\rm By construction, the presented algorithm produces upon termination (if any)
\begin{itemize}
\item an execution protocol $y^t=\{y_\tau,g(y_\tau)\}_{\tau=1}^t$, where $t$ is the step where the algorithm terminates, and $y_\tau$, $1\leq \tau\leq t$, are the search points generated in course of the run; by construction, all these search points belong to $Y$;
\item an accuracy certificate $\lambda^t$ -- a collection of nonnegative weights $\lambda_1,...,\lambda_t$ summing up to 1 -- such that the affine function
$
h(y)=\sum_{\tau=1}^t \lambda_\tau \langle g(y_\tau),y_\tau-y\rangle
$
satisfies the relation
$
\epsilon(y^t,\lambda^t):=\max\limits_{x\in Y}h(x)\leq\epsilon,
$
where $\epsilon$ is the target tolerance, exactly as required in (\ref{bruce1}).
\end{itemize}
}
\end{remark}
\subsubsection{Efficiency estimate}
\begin{proposition}\label{PropNERML}
Given on input a target tolerance $\epsilon>0$, the NERML algorithm terminates after finitely many steps, with execution protocol $y^t$ and accuracy certificate $\lambda^t$, described in Remark \ref{remB}. The number of steps of the algorithm does not exceed
\begin{equation}\label{numberofsteps}
N=C(\gamma,\theta){\Omega^2L^2[g]\over\epsilon^2} \hbox{, where } C(\gamma,\theta)={(1+\gamma^2)\over\gamma^2[1-[\gamma+(1-\gamma)\theta]^2]}.
\end{equation}
\end{proposition}
For proof, see Section \ref{Appendix.PropNERML}.
\begin{remark}\label{remNERML2} {\rm Inspecting the proof of Proposition \ref{PropNERML}, it is immediately seen that when $\omega(\cdot)$ is continuously differentiable on the entire $Y$, one can replace the rule (\ref{specify}) with
$$
u_{\tau+1}=\argmin_{y\in Y_\tau} [\omega(y)-\langle \omega'(y^s),y-y^s\rangle],
$$
where $y^s$ is an arbitrary point of $Y^o=Y$. The cost of this modification is that of replacing $\Omega$ in the efficiency estimate with $\Omega^+$, see Remark \ref{remLevel}. Computational experience shows that a good choice of $y^s$ is the best, in terms of the objective, search point generated before the beginning of phase $s$.}
\end{remark}
\paragraph{Solving $(P)$ and $(D)$ by NERML} is completely similar to the case of MDL, with the bound
\begin{equation}\label{NERMLbound}
t(\epsilon)\leq \hbox{\rm Ceil}\left({\Omega^2L_f^2(1+\gamma^2)\over\gamma^2[1-[\gamma+(1-\gamma)\theta]^2\epsilon^2}\right)\quad [\Omega=\Omega[Y,\omega(\cdot)]]
\end{equation}
in the role of (\ref{MDLbound}).
\begin{remark}\label{reminaccurate} {\rm Observe that Propositions \ref{PropMD}-\ref{PropNERML} do not impose restrictions of the vector field $g(\cdot)$ processed by the respective algorithms aside from the boundedness assumption (\ref{normalization}). Invoking Remark \ref{rem1}, we arrive at the following conclusion:
\par
{\sl in the situation of section \ref{sect:situation} and given $\delta\geq0$,  let instead of exact maximizers $x(y)\in\Argmax_{x\in X}\langle x,Ay+a\rangle$, approximate maximizers $x_\delta(y)\in X$ such that $\langle x_\delta(y),Ay+a\rangle \geq \langle x(y),Ay+a\rangle-\delta$ for all $y\in Y$ be available. Let also
$$
f^\prime_\delta(y)=A^*x_\delta(y)+\psi'(y)
$$
be the associated approximate subgradients of the objective $f$ of $(D)$. Assuming
$$
L_{f,\delta} =\sup\limits_{y\in Y} \|f^\prime_\delta(y)\|<\infty,
$$
let MD/MDL/NERML be applied to the vector field $g(\cdot)=f^\prime_\delta(\cdot)$. Then the number of steps of each method before termination  remains bounded by the respective bound {\rm (\ref{MDini}), (\ref{numberofstepsLevel}) or (\ref{numberofsteps})}, with $L_{f,\delta}$ in the role of $L_f$. Besides this, defining the approximate solutions to $(P)$, $(D)$ according to {\rm (\ref{eq229})}, with $x_\tau=x_\delta(y_\tau)$, we ensure the validity of $\delta$-relaxed version of
the accuracy guarantee {\rm (\ref{accuracy})}, specifically, the relation
$$
[f(\widehat{y}^t)-\Opt(P)]+[\Opt(D)-f_*(\widehat{x}^t)]\leq \epsilon+\delta.
$$
}}
\end{remark}
\section{An alternative: Smoothing}\label{sect:smoothing}\label{sect:SCG}
An alternative to the approach we have presented so far is based on the use of the proximal setup for $Y$ to {\sl smoothen} $f_*$ and then to maximize the resulting smooth approximation of $f_*$  by the Conditional Gradient (CG) algorithm. This approach is completely similar to the one used by Nesterov in his breakthrough paper \cite{NesSmooth}, with the only difference that since in our situation domain $X$ admits LO oracle rather than a good proximal setup, we are bounded to replace the $O(1/t^2)$-converging Nesterov's method for smooth convex minimization with $O(1/t)$-converging CG.

\par
Let us describe the CG implementation in our setting. Suppose that we are given a norm $\|\cdot\|_x$ on $E_x$, a representation (\ref{Fenchel}) of $f_*$, a proximal point setup $(\|\cdot\|_y,\;\omega(\cdot))$ for $Y$ and a desired tolerance $\epsilon>0$. We assume w.l.o.g. that $\min\limits_{y\in Y} \omega(y)=0$ and set, following Nesterov \cite{NesSmooth},
\begin{equation}\label{apprf}
\begin{array}{c}
f_*^\beta(x)=\min\limits_{y\in Y}\left[\langle x,Ay+a\rangle + \psi(y)+\beta\omega(y)\right] \\
\beta=\beta(\epsilon):={\epsilon\over \Omega^2},\,\, \Omega=\Omega[Y,\omega(\cdot)]\\
\end{array}
\end{equation}
From (\ref{Fenchel}), the definition of $\Omega$ and the relation $\min_Y\omega=0$ it immediately follows that
\begin{equation}\label{close}
\forall x\in X: f_*(x)\leq f_*^\beta(x)\le f_*(x)+{\epsilon\over 2},
\end{equation}
and $f_*^\beta$ clearly is concave. It is well known (the proof goes back to J.-J. Moreau \cite{Mor62,Mor65}) \label{Moreau2} that strong convexity, modulus 1 w.r.t. $\|\cdot\|_y$, of $\omega(y)$ implies smoothness of $f_*^\beta$, specifically,
\begin{equation}\label{smoothness}
\forall (x,x'\in X): \;\|\nabla f_*^\beta(x)-\nabla f_*^\beta (x')\|_{x,*}\leq {1\over\beta}\|A\|_{y;x,*}^2\|x-x'\|_x,
\end{equation}
where
\[
\|A\|_{y;x,*}=\max\{\|A^*u\|_{y,*}:u\in E_x,\|u\|_x\leq1\}=\max\{\|Ay\|_{x,*}:y\in E_y,\|y\|_y\leq1\}.\]
Observe also that under the assumption that an optimal solution $y(x)$ of the right hand side minimization problem in (\ref{apprf}) is available at a moderate computational cost, \footnote{In typical applications, $\psi$ is just linear, so that computing $y(x)$ is as easy as computing the value of the prox-mapping associated with   $Y$, $\omega(\cdot)$.} we have at our disposal a FO oracle for $f_*^\beta$:
$$
f_*^\beta (x)=\langle x,Ay(x)+a\rangle+\psi(y(x))+\beta\omega(y(x)),\,\,\ \nabla f_*^\beta(x)=Ay(x)+a.
$$
We can now use this oracle, along with the LO oracle for $X$, to solve $(P)$ by CG\footnote{ If our only goal were to approximate $f_*$ by a smooth concave function, we could use other options, most notably the famous Moreau-Yosida regularization \cite{Mor65,Yos64} \label{Moreau3}. The advantage of the smoothing based on Fenchel-type representation  (\ref{Fenchel}) of $f_*$ and proximal setup for $Y$ is that in many cases it indeed allows for computationally cheap approximations, which is usually not the case  for Moreau-Yosida regularization.}. In the sequel, we refer to the outlined algorithm as to SCG (Smoothed Conditional Gradient).
\paragraph{Efficiency estimate} for SCG is readily given by Proposition \ref{propredgrad}. Indeed, assume from now on that $X$ is contained in $\|\cdot\|_x$-ball of radius $R$ of $E_x$. It is immediately seen that under this assumption, (\ref{smoothness}) implies the validity of the condition (cf. (\ref{suchthat}) with $q=2$)
\begin{equation}\label{conditionagain}
\begin{array}{c}
\forall x,x'\in X: f_*^\beta(x')\geq f_*^\beta(x)+\langle\nabla f_*^\beta(x),x'-x\rangle -{1\over 2}\cL\|x'-x\|_X^2,\quad
\cL={R^2\|A\|_{y;x,*}^2\over\beta}={R^2\|A\|_{y;x,*}^2\Omega^2\over \epsilon}.
\end{array}
\end{equation}
In order to find an $\epsilon$-maximizer of $f_*$, it suffices, by (\ref{close}), to find an $\epsilon/2$-maximizer of $f_*^\beta$; by (\ref{then}) (where one should set $q=2$), what takes
\begin{equation}\label{redgrad}
t^{\cg}(\epsilon) =O(1)\epsilon^{-1}\cL=O(1){\Omega^2\|A\|_{y;x,*}^2R^2\over\epsilon^2}
\end{equation}
steps.
\paragraph{Discussion.} Let us assume, as above, that $X$ is contained in the centered at the origin $\|\cdot\|_x$-ball of radius $R$, and let us compare the essentially identical to each other\footnote{provided the parameters $\gamma\in(0,1)$, $\theta\in (0,1)$ in (\ref{MDLbound}), (\ref{NERMLbound}) are treated as absolute constants.} complexity bounds (\ref{MDefficiency}), (\ref{MDLbound}), (\ref{NERMLbound}), with the bound (\ref{redgrad}).  Under the natural assumption that the subgradients of $\psi$ we use satisfy the bounds $\|\psi'(y)\|_{y,*}\leq L_\psi$, where $L_\psi$ is the Lipschitz constant of $\psi$ w.r.t. the norm $\|\cdot\|_y$, (\ref{Lf}) implies that
\begin{equation}\label{Lfagain}
L_f\leq \|A\|_{y;x,*}R+L_\psi.
\end{equation}
Thus, the first three complexity bounds reduce to
\begin{equation}\label{MDbounds}
{\cal C}_{\md}(\epsilon)=\hbox{\rm Ceil}\left(O(1){[R\|A\|_{y;x,*}+L_\psi]^2\Omega^2[Y,\omega(\cdot)]\over\epsilon^2}\right),
\end{equation}
while the conditional gradients based complexity bound is
\begin{equation}\label{smoothbounds}
{\cal C}_{\cg}(\epsilon)=\hbox{\rm Ceil}\left(O(1){[R\|A\|_{y;x,*}]^2\Omega^2[Y,\omega(\cdot)]\over\epsilon^2}\right).
\end{equation}
We see that {\sl assuming $L_\psi\leq O(1)R\|A\|_{y;x,*}$} (which indeed is the case in many applications, in particular, in the examples we are about to consider), {\sl the complexity bounds in question are essentially identical.} This being said, we believe that the two approaches in question seem to have their own advantages and disadvantages. Let us name just a few:
\begin{itemize}
\item Formally, the SCG has a more restricted area of
applications than MD/MDL/NERML, since relative simplicity  of the optimization problem  in  (\ref{apprf}) is a more restrictive requirement than relative simplicity of computing prox-mapping associated with $Y,\omega(\cdot)$. At the same time, in most important applications known to us $\psi$ is just linear, and in this case the just outlined phenomenon disappears.
\item An argument in favor of SCG is its insensitivity to the Lipschitz constant of $\psi$. Note, however, that in the case of linear $\psi$ (which, as we have mentioned, is the case of primary interest) the nonsmooth techniques admit simple modifications (not to be considered here) which make them equally insensitive to $L_\psi$.
\item Our experience shows that the convergence pattern of nonsmooth methods utilizing memory (MDL and NERML) is, at least at the beginning of the solution process,  much better than is predicted by their worst-case efficiency estimates.
    It should be added that in theory there exist situations where the nonsmooth approach ``most probably,''  or even provably, significantly outperforms the smooth one. This is the case when $E_y$ is of moderate dimension. A well-established {\sl experimental} fact is that when solving $(D)$ by MDL, every $\dim E_y$ iterations of the method reduce the inaccuracy by an absolute constant factor, something like 3. It follows that if $n$ is in the range of few hundreds, a couple of thousands of MDL steps can yield a solution of accuracy which is incomparably better than the one predicted by the theoretical worst-case oriented  $O(1/\epsilon^2)$ complexity bound of the algorithm. Moreover, in principle one can solve $(D)$ by the Ellipsoid method  with certificates \cite{NOR}, building accuracy certificate of resolution $\epsilon$ in {\sl polynomial time} $O(1)n^2\ln(L_f\Omega[Y,\omega(\cdot)]/\epsilon)$. It follows that when $\dim E_y$ is in the range of few tens, the nonsmooth approach allows to solve, in moderate time, problems  $(P)$ and $(D)$ to high accuracy.
    Note that low dimensionality of $E_y$ by itself does not prevent $X$ to be high-dimensional and ``difficult;'' how frequent are these situations in actual applications, this is another story.
\end{itemize}
We believe that the choice of one, if any, of the outlined approaches to use, is the issue which should be resolved, on the case-by-case basis, by computational practice. We believe, however, that it makes sense to keep them both in mind.
\section{Application examples}\label{sect:examples}
In this section we work out some application examples, with the goal to demonstrate that the approach we are proposing
possesses certain application potential.
\subsection{Uniform norm matrix completion}
Our first example (for its statistical motivation, see \cite{Fan11}) is as follows: given a symmetric $p\times p$ matrix $b$ and a positive real $R$, we want to find the best entrywise approximation of $b$
by a positive semidefinite matrix $x$ of given trace $R$, that is, to solve the problem
\begin{equation}\label{lastexample}
\begin{array}{c}
\min_{x\in X} \left[-f_*(x)=\|x-b\|_\infty\right]\\
X=\{x\in\bS^p:\; x\succeq0,\;\Tr(x)=R\},\;\|x\|_\infty=\max\limits_{1\leq i, j\leq p} |x_{ij}|\\
\end{array}
\end{equation}
where $\bS^p$ is the space of $p\times p$ symmetric matrices.
    Note that with our $X$,  computing prox-mappings associated with all known proximal setups needs eigenvalue decomposition of a $p\times p$ symmetric matrix and thus becomes
    computationally demanding in the large scale case. On the other hand, to maximize a linear form $\langle \xi,x\rangle=\Tr(\xi x)$ over $x\in X$ requires computing the maximal eigenvalue of $\xi$ along with corresponding eigenvector. In the large scale case this task is by orders of magnitude less demanding than computing full eigenvalue decomposition. Note that our $f_*$ admits a simple Fenchel-type representation:
$$
f_*(x)=-\|x-b\|_\infty=\min\limits_{y\in Y} \left[f(y)=\langle -x,y\rangle +\langle b,y\rangle\right], \,\,Y=\{y\in\bS^p: \|y\|_1:=\sum_{i,j}|y_{ij}|\leq1\}.
$$
Equipping $E_y=\bS^p$ with the norm $\|\cdot\|_y=\|\cdot\|_1$, and $Y$  with the d.-g.f.
$$
\omega(y)=\alpha\ln(p)\sum_{i,j=1}^p|y_{ij}|^{1+r(p)},\,\,r(p)={1\over\ln(p)},
$$
where $\alpha$ is an appropriately chosen constant of order of 1 (induced by the necessity to make $\omega(\cdot)$ strongly convex, modulus 1, w.r.t. $\|\cdot\|_1$), we
get a proximal setup for $Y$ such that
$$
\Omega[Y,\omega(\cdot)]\leq O(1)\sqrt{\ln p}.
$$
We see that our problem of interest fits well the setup of methods developed in this paper. Invoking the bounds (\ref{MDbounds}), (\ref{smoothbounds}), we conclude that (\ref{lastexample}) can be solved within accuracy $\epsilon$ in at most
\begin{equation}\label{atmost65}
t(\epsilon)=O(1){[R+\|b\|_\infty]^2\ln(p)\over \epsilon^2}
\end{equation} steps by any of methods MD, MDL or NERML, and in at most $O(1){R^2\ln(p)\over \epsilon^2}$ steps by SCG.
\par
It is worth to mention that in the case in question, the algorithms yielded by the nonsmooth approach admit a ``sparsification'' as follows.  We are in the case of $\langle x,Ay+a\rangle \equiv  \Tr(xy)$, and $X=\{x:x\succeq0, \Tr(x)=R\}$, so that $x(y)=R e_ye_y^T$, where $e_y$ is the leading eigenvector of a matrix $y$ normalized to have $\|e_y\|_2=1$. Given a desired accuracy $\epsilon>0$ and a unit vector $e_{y,\epsilon}$ such that $\Tr(y[e_{y,\epsilon}e_{y,\epsilon}^T])\geq \Tr(y[e_ye_y^T])-R^{-1}\epsilon$, and setting $x_\epsilon(y)=Re_{y,\epsilon}e_{y,\epsilon}^T$, we ensure that $x_\epsilon(y)\in X$ and that $x_\epsilon(y)$ is an $\epsilon$-maximizer of $\langle x,Ay+a\rangle$ over $x\in X$. Invoking
Remark \ref{reminaccurate}, we conclude that when utilizing $x_\epsilon(\cdot)$ in the role of $x(\cdot)$, we get $2\epsilon$-accurate solutions to $(P)$, $(D)$ in no more than $t(\epsilon)$ steps. Now, we can take as $e_{y,\epsilon}$ the normalized leading eigenvector of an arbitrary matrix $\widehat{y}(y)$ such that $\|\sigma(\widehat{y}-y)\|_{\infty}\leq\epsilon$. Assuming $R/\epsilon>1$ and given $y\in Y$, let us sort the magnitudes of entries in $y$ and build $y_\epsilon$ by ``thresholding'' -- by zeroing out as many smallest in magnitude entries as possible under the restriction that the remaining part of the matrix $y$ is symmetric, and the sum of squares of the entries we have replaced with zeros does not exceed $R^{-2}\epsilon^2$. Since $\|y\|_1\leq1$, the number $N_\epsilon$ of nonzero entries in $y_\epsilon$ is at most $O(1)R^2/\epsilon^2$. On the other hand,
by construction, the Frobenius norm $\|\sigma(y-y_\epsilon)\|_2$ of $y-y_\epsilon$ is $\le R^{-1}\epsilon$, thus
$\|\sigma(y-y_\epsilon)\|_{\infty}\leq R^{-1}\epsilon$, and we can take as $e_{y,\epsilon}$ the normalized leading eigenvector of $y_\epsilon$. When the size $p$ of $y$ is $\gg R^2/\epsilon^2$ (otherwise the outlined sparsification does not make sense), this approach reduces the problem of computing the leading eigenvector to the case when the matrix is question is relatively sparse, thus reducing its computational cost.
\subsection{Nuclear norm SVM}
Our next example is as follows: we are given an $N$-element  sample of $p\times q$ matrices $z_j$ (``images'')  equipped with labels $\epsilon_j\in\{-1,1\}$. We assume the images to be normalized by the restriction \begin{equation}\label{mu}
\|\sigma(z_j)\|_\infty\leq1.
\end{equation}
 We want to find a linear classifier of the form
$$
\widehat{\epsilon}_j=\sign(\langle x, z_j\rangle + b).\eqno{[\langle z,x\rangle=\Tr(zx^T)]}
$$
which predicts well labels of images. We assume that there exist such right and left orthogonal transformations of the image, that the label can be predicted using only a small number of diagonal elements of the transformed image. This implies that the classifier we are looking for is sparse in the corresponding basis, or that the matrix $x$ is of low rank. We arrive at the ``low-rank-oriented'' SVM-based reformulation of this problem:
\begin{equation}\label{ini}
\min\limits_{x:\|\sigma(x)\|_1\leq R}\left[h(x):=\min_{b\in \bR}\left[N^{-1} {\sum}_{j=1}^N\left[1 -\epsilon_j\left[\langle x,z_j\rangle +b\right]\right]_+\right]\right],
\end{equation}
where $\|\sigma(\cdot)\|_1$ is the nuclear norm, $[a]_+=\max[a,0]$, and $R\geq1$ is a parameter.\footnote{The restriction $R\geq1$ is quite natural. Indeed, with the optimal choice of $x$, we want most of the terms $\left[1 -\epsilon_j\left[\langle x,z_j\rangle +b\right]\right]_+$ to be $\ll 1$; assuming that the number of examples with $\epsilon_j=-1$ and $\epsilon_j=1$ are of order of $N$, this condition can be met only when $|\langle x,z_j\rangle|$ are at least of order of 1 for most of $j$'s. The latter, in view of (\ref{mu}), implies that $\|\sigma(x)\|_{1}$ should be at least $O(1)$.}
\par
In this case the domain $X$ of problem ($P$) is the ball of the nuclear norm in the space $\bR^{p\times q}$ of $p\times q$ matrices and $p,q$ are large. As we have explained in the introduction, same as in the example of the previous section, in this case the computational complexity of LO oracle is typically much smaller than the complexity of computing prox-mapping. Thus, from practical viewpoint, in a meaningful range of values of $p,q$ the LO oracle is ``affordable,''  while the prox-mapping is not.
\par
Observing that $[a]_+=\max_{0\leq y\leq 1}ya$, and denoting  $\mones=[1;...;1]\in
\bR^q$, we get
$$
N^{-1} \sum_{j=1}^N\left[1 -\epsilon_j\left[\langle x,z_j\rangle +b\right]\right]_+=
\max\limits_{y:0\leq y\leq \mones}\left\{N^{-1}\sum_{j=1}^Ny_j\left[1-\epsilon_j\left[\langle x,z_j\rangle +b\right]\right]\right\},
$$
whence
\begin{equation}\label{h}
h(x)=\max\limits_{y\in Y}N^{-1}\sum_{j=1}^Ny_j\left[1-\epsilon_j\langle x,z_j\rangle\right],
\end{equation}
where
\begin{equation}\label{Y}
Y=\{y\in\bR^N: 0\leq y\leq \mones, \sum_j\epsilon_jy_j=0\};
\end{equation}
from now on we assume that $Y\neq\emptyset$. When setting
\begin{equation}\label{setting}
\cA y =N^{-1}\sum_{j=1}^Ny_j\epsilon_jz_j: \bR^N\to\bR^{p\times q},\,\, X=\{x\in\bR^{p\times q}: \|\sigma(x)\|_{1}\leq R\}, \,\,\psi(y)=-N^{-1} \mones^Ty
\end{equation}
and passing from minimizing $h(x)$ to maximizing $f_*(x)\equiv -h(x)$, problem (\ref{ini}) becomes
\begin{equation}\label{inireformulated}
\max\limits_{x\in X}\left[f_*(x):=\min\limits_{y\in Y}\left[\langle x,\cA y\rangle + \psi(y)\right]\right].
\end{equation}
Let us equip $E_y=\bR^N$ with the standard Euclidean norm $\|\cdot\|_2$, and $Y$ - with the Euclidean d.-g.f. $\omega(y)={1\over 2}y^Ty$. Observe that
$$
\left[f^\prime(y)\right]_j=N^{-1}\left[\langle x(y),\epsilon_jz_j\rangle -1\right],\;\; x(y)\in\Argmax_{x\in X} \Tr(xy^T),
$$
meaning that
$$
\|f^\prime(y)\|_\infty\leq \max_j\left[N^{-1}\left[1+\|\sigma(x(y))\|_{1}\|\sigma(z_j)\|_{\infty}\right]\right]\leq N^{-1}[R+1]\leq 2N^{-1}R.
$$
Using our notation of section \ref{sect:certificates} we have
$$
L_f:=\sup_{y\in Y} \|f^\prime(y)\|_*\leq 2N^{-1/2}R
$$
(we are in the case of $\|\cdot\|_*=\|\cdot\|_2$), and,
besides,
$$
\Omega_Y\leq \sqrt{N/2}.
$$
 We conclude that for every $\epsilon>0$, the number $t$ of MD steps needed to ensure (\ref{accuracy}) does not exceed
 $$
 t_{\md}(\epsilon)=\hbox{\rm Ceil}\left({2R^2\over\epsilon^2}\right)
 $$
 (see (\ref{MDefficiency})), and similarly for MDL, NERML, and SCG.

\subsection{Multi-class classification under $\infty|2$ norm constraint}
Our last example illustrates the potential of the proposed approach in the case when the domain $X$ of $(P)$ does not admit a proximal setup with ``moderate'' $\Omega_X$. Namely,  let
$(P)$ be the problem
\begin{equation}\label{Pnowis}
\min_{x\in X} \left[-f_*(x):=\|Bx -b\|_{y,*}\right] \qquad [x\mapsto Bx: E_x\to E_y]
\end{equation}
where $\|\cdot\|_{y,*}$ is the norm conjugate to a norm $\|\cdot\|_y$ on $E_y$. We are interested in the case of box-type $X$, specifically,
\begin{equation}\label{Xis}
X=\{[x^1;...;x^M]\in\bR^{n_1}\times...\times\bR^{n_M}: \|x^i\|_2\leq R,\,1\leq i\leq M\}
\end{equation}
As it was mentioned in Introduction, for every proximal setup $(\|\cdot\|,\omega_x(\cdot))$ for $X$ which is normalized by the requirement that simple ``well behaved'' on $X$ convex functions should have moderate Lipschitz constants w.r.t. $\|\cdot\|$ (specifically, the coordinates of $x\in X$ should have Lipschitz constants $\leq1$), one has  $\Omega[X,\omega_x(\cdot)]\geq O(1)\sqrt{M}R$. As a result, the theoretical complexity of the FO methods as applied to (\ref{Pnowis}) grows with $M$ at the rate at least $O(\sqrt{M})$, thus becoming prohibitively high for large $M$. We are about to show that the approaches developed in this paper are free of this shortcoming. Specifically, we can easily build a Fenchel-type representation of $f_*$:
$$
f_*(x)=-\|Bx-b\|_{y,*}=\min_{y\in Y}\left[\langle B^*y,x\rangle -\langle b,y\rangle\right], \, Y=\{y\in E_y:\|y\|_y\leq1\}.
$$
Assume that $Y$ admits a good proximal setup.  We can
augment $\|\cdot\|_y$ with a d.-g.f. $\omega_y(\cdot)$ for $Y$ such that $\|\cdot\|_y,\omega(y)$ form a proximal setup, and applying any of the methods we have developed in sections \ref{sect:certificates} and \ref{sect:smoothing}, the complexity of finding $\epsilon$-solution to (\ref{Pnowis}) by any of these methods becomes
\[
O(1)\left({R\|B\|_{x;y,*}+\|b\|_{y,*}\over\epsilon}\right)^2, \,\,\|B\|_{x;y,*}=\max_{x=[x^1;...;x^M]}\left\{\|Bx\|_{y,*}: \|x\|_{\infty|2}:=\max_i\|x^i\|_2\leq1\right\}.
\]
Note that in this bound $M$ does not appear, at least explicitly.
\paragraph{Multi-class classification problem} we consider is as follows:  we observe $N$ ``feature vectors'' $z_j\in\bR^q$, each belonging to one of $M$ non-overlapping classes, along with labels $\chi_j\in\bR^M$ which are basic orths in $\bR^M$; the index of the (only) nonzero entry in $\chi_j$ is the number of class to which $z_j$ belongs. We want to build a multi-class analogy of the standard linear classifier as follows: a multi-class classifier is specified by a matrix $x\in\bR^{M\times q}$ and a vector $b\in\bR^M$. Given a feature vector $z$, we compute the $M$-dimensional vector $x z+b$, identify its maximal component, and treat the index of this component as our guess for the serial number of the class to which $z$ belongs.\par
The multi-class analogy of the usual approach to building binary classifiers by minimizing the empirical hinge loss is as follows \cite{Cramer01,Amit:etal:2007}. Let $\bar{\chi}_j=\mones-\chi_j$ be the ``complement'' of $\chi_j$.
Given a feature vector $z$ and the corresponding label $\chi$, let us set
$$
h=h(x,b;z,\chi)=[x z+b]-[\chi^T[x z +b]]\mones +\bar{\chi}\in \bR^M\eqno{[\mones=[1;...;1]\in\bR^M]}.
$$
Note that if $i_*$ is the index of the only nonzero entry in $\chi$, then the $i_*$-th entry in $h$ is zero (since $\chi_{i_*}=1$). Further, $h$ is nonpositive
 if and only if the classifier, given by $x,b$ and evaluated at $z$, ``recovers the class $i_*$ of $z$ with margin 1'', i.e., we have $[xz+b]_j\le [xz+b]_{i_*}-1$  for  $j\neq i_*$.  On the other hand, if the classifier fails to classify $z$ correctly (that is, $[xz+b]_j \geq [xz+b]_{i_*}$ for some $j\neq i_*$), then the maximal entry in $h$ is $\geq 1$. Altogether, when  setting
$$
\eta(x,b;z,\chi)=\max_{1\leq j\leq M} [h(x,b;z,\chi)]_j,
$$
we get a nonnegative function which vanishes for the pairs $(z,\chi)$ which are ``quite reliably'' -- with margin $\geq1$ -- classified by $(x,b)$, and is $\geq1$ for the pairs $(z,\chi)$ with $z$ not classified correctly. Thus the function
$$
H(x,b)=\bE\{\eta(x,b;z,\chi)\},
$$
the expectation being taken over the distribution of examples $(z,\chi)$, is an upper bound on the probability for classifier $(x,b)$ to misclassify a feature vector. What we would like to do now is to minimize $H(x,b)$ over $x,b$. To do this, since $H(\cdot)$ is not observable, we replace the expectation by its empirical counterpart
$$
H_N(x,b)=N^{-1}\sum_{j=1}^N\eta(x,b;z_j,\chi_j).
$$
For the sake of simplicity (and, upon a close inspection, without much harm), we assume from now on that $b=0$.\footnote{To arrive at this situation, one can augment  $z_j$ by additional entry, equal to 1, and to redefine $x$: the new $x$ is the old  $[x,b]$.} Imposing, as it is always the case in hinge loss optimization, an upper bound on some norm $\|x\|_x$ of $x$, we arrive at the optimization problem
\begin{equation}\label{multiclass}
\min_{x\in X} \left[-f_*(x)=N^{-1}\sum_{j=1}^N \max\limits_{i\leq M} [x z_j-[\chi_j^Tx z_j]\mones+\bar{\chi}_j]_i\right], \quad X=\{x:\|x\|_x\leq R\}.
\end{equation}
From now on we assume that $z_j$'s are normalized:
\begin{equation}\label{zj}
\|z_j\|_2\leq1,\,1\leq j\leq N.
\end{equation}
Under this constraint, a natural (although not the only meaningful) choice of the norm $\|\cdot\|_x$ is the maximum of the $\|\cdot\|_2$-norms of the rows $[x^i]^T$ of $x$. If we identify $x$ with the vector $[x^1;...;x^M]$, $X$ becomes the set (\ref{Xis}) with $n_1=n_2=...=n_M=q$, and the norm $\|\cdot\|_x$ becomes $\|\cdot\|_{\infty|2}$. The same argument as in the previous section allows us to assume that $R\geq1$.
\par
Noting that$\min\limits_ih_i=\min_u\{u^Th: u\geq0,\sum_iu_i=1\}$, (\ref{multiclass}) can be rewritten as
\begin{equation}\label{multiclassbecomes}
\max\limits_{x\in X}\left[f_*(x)=\min\limits_{y\in Y} [\langle y,Bx\rangle +\psi(y)]\right]
\end{equation}
where
$$
\begin{array}{l}
Y=\{y=[y^1;...;y^N]:\, y^j\in\bR^M_+,\,\sum_i[y^j]_i=N^{-1},\,1\leq j\leq N\} \subset E_y=\bR^{MN}\\
Bx=[B^1x;...;B^Nx],\;B^jx=\left[z_j^T[x^{i(j)}-x^1];...;z_j^T[x^{i(j)}-x^M]\right],\,j=1,...,N\\
\psi(y)=\psi(y^1,...,y^N)=-\sum_{j=1}^N[y^j]^T\bar{\chi}_j;\\
\end{array}
$$
(here $i(j)$ is the class of $z_j$, i.e., the index of the only nonzero entry in $\chi_j$). Note that $Y$ is a part of the standard simplex $\Delta_{MN}=\{y\in\bR^{MN}_+:\,\sum_{j=1}^N\sum_{i=1}^{M} [y^j]_i=1\}\subset E_y=\bR^{MN}$. Equipping $E_y$ with the norm $\|\cdot\|_y=\|\cdot\|_1$ (so that $\|\cdot\|_{y,*}=\|\cdot\|_\infty$), and $Y$ -- with the entropy d.-g.f.
$$
\omega_y(y)=\sum_{j=1}^N\sum_{i=1}^M [y^j]_i\ln([y^j]_i)
$$
(known to complete $\|\cdot\|_1$ to a proximal setup for $\Delta_{MN}$),
we get a proximal setup for $Y$ with $\Omega=\Omega[Y,\omega(\cdot)]\leq \sqrt{2\ln(M)}$. Next, assuming $\|x\|_x\equiv \|x\|_{\infty|2}\leq1$, we have
$$
\|Bx\|_{y,*}=\|Bx\|_\infty=\max\limits_{{1\leq i\leq m\atop1\leq j\leq N}}|z_j^T[x^{i(j)}-x^i]|\leq \|z_j\|_2\|x^{i(j)}-x^i\|_2\leq 2,
$$
so that $\|B\|_{x;y,*}\leq2$. Furthermore, $\psi$ clearly is Lipschitz continuous with constant 1 w.r.t. $\|\cdot\|_y=\|\cdot\|_1$. It follows that the complexity of finding an $\epsilon$-solution to (\ref{multiclassbecomes}) by MD, MDL, NERML or SCG is bounded by $O(1){R^2\ln(M)\over\epsilon^2}$ (see (\ref{Lfagain}), (\ref{MDbounds}), (\ref{smoothbounds}) and take into account that $R\geq1$, and that what is now called $B$, was called  $A^*$ in the notation used in those bounds, so that $\|B\|_{x;y,*}=\|A\|_{y;x,*}$). Note that the resulting complexity bound is independent of $N$ and is ``nearly independent'' of $M$. Finally, prox-mapping for $Y$ is given by a closed form expression and can be computed in linear time:
$$
\begin{array}{l}
\argmin\limits_{\{y^j\in\bR^M\}_{j=1}^N}\left\{\sum_{j=1}^N\sum_{i=1}^M[y^j]_i\ln([y^j]_i)+\sum_{j=1}^N\langle \xi^j,y^j\rangle:y^j
\geq0,\sum_{i=1}^M[y^j]_i=N^{-1},\;1\leq j\leq N\right\}\\
=\left\{\widehat{y}^j: [\widehat{y}^j]_i={\exp\{-[\xi^j]_i\}\over N\sum_{s=1}^M\exp\{-[\xi^j]_s\}},\,1\leq i\leq M\right\}_{j=1}^N.\\
\end{array}
$$
\section{NERML: Numerical illustration}\label{sect:numerics}
The goal of the numerical experiments to be reported is to illustrate how the performance of NERML scales up with the dimensions  of $x$ and $y$ and the memory of the method. Below we consider a kind of matrix completion problem, specifically,
\be
\Opt=\min\limits_{x\in\bR^{p\times p}}\left\{\|\cP(x-a)\|_\infty:\;\|\sigma(x)\|_1\leq1\right\}.
\ee{pprim_num}
Here $a$ is a given matrix, $\cP$ is a linear mapping from $\bR^{p\times p}$ into $\bR^N$, and $\|y\|_\infty=\max_i|y_i|$ is the uniform norm on $\bR^N$. In our experiments, $\cP$  is defined as follows. We select a set $\cI=\{(i,j)\}$ of $rp$ cells in a $p\times p$ matrix in such a way that every row and every column contains exactly $r$ of the selected cells. We then label at random the selected cells by indexes from $\{1,2,...,N\}$, with the only restriction that every one of the $N$ indexes labels the same number $pr/N$ (which with our choice of $r,p,N$ always is integer) of the cells. The $i$-th, $1\leq i\leq N$, entry in $\cP y$ is the sum, over all cells from $\cI$ labeled by $i$, of the entries of $y\in\bR^{p\times p}$ in the cells. With $N=pr$, $\cP y$ is just the restriction of $y$ onto the cells from $\cI$, and (\ref{pprim_num}) is the standard matrix completion problem with uniform fit. Whatever simplistic and ``academic,'' our setup, in accordance with the goals of our numerical study, allows for full control of image dimension of $\cP$ (i.e., the design dimension of the problem actually solved by NERML) whatever large be the matrices $x$ and the set $\cI$.
\par
Our test instances were generated as follows: given $p,r,N$, we generate at random the set $\cI$ along with its labeling (thus specifying $\cP$) and a vector $w\in \bR^N$ with $d=32$ nonzero entries. Finally, we set
$$
v=\|\sigma(\cP^*w)\|_1^{-1}\cP^*w,\;\;\mbox{and}\;\;a=v+2\|v\|_\infty\xi,\,\,\,\|v\|_\infty=\max_{1\leq i,j\leq p}|v_{ij}|,
$$
where the entries in $p\times p$ ``noise matrix'' $\xi$ are the projections onto $[-1,1]$ of random reals sampled, independently of each other, from the standard Gaussian distribution.
\par
Written in the form of $(P)$, problem (\ref{pprim_num}) reads
\begin{equation}\label{pprim_num_P}
[-\Opt]\ = \max_{x\in X} \left\{f_*(x)=-\|\cP(x-a)\|_\infty\right\},\,\,X=\{x\in\bR^{p\times p}:\|\sigma(x)\|_1\leq1\};
\end{equation}
the Fenchel-type representation (\ref{Fenchel}) of $f_*$ is
$$
f_*(x)=\min\limits_{{y\in\bR^N,\atop\|y\|_1\leq1}} y^T\cP(x-a),
$$
so that the dual problem to be solved by NERML is
\begin{equation}\label{dual_num}
\min_{y\in Y} \left\{f(y)=\max\limits_{x:\|\sigma(x)\|_1\leq1}\langle x,\cP^*y\rangle -[\cP a]^Ty\right\}, \,\, Y=\{y\in\bR^N: \|y\|_1\leq1\}.
\end{equation}
\par
Before passing to numerical results, we make the following important remark. Our description of NERML in section \ref{NERML} is adjusted to the case when the accuracy $\epsilon$ to which (\ref{pprim_num_P}) should be solved is given in advance.
Note, however, that $\epsilon$ is used {\sl only} in the termination rule {\bf B.2.2.3}: we terminate when the optimal value in the current auxiliary problem (\ref{bruce10}) becomes $\leq\epsilon$. The optimal value in question is a certain ``online observable'' function  $\epsilon_\tau$ of the ``time'' $\tau$ defined as the total number of steps performed so far. Moreover, at every time $\tau$ we have at our disposal a feasible solution $x_\tau$ to the problem of interest $(P)$ (in our current situation, to (\ref{pprim_num_P})) such that
$$
\Opt-f_*(x_\tau)\leq\epsilon_\tau
$$
-- this is the solution which NERML, as defined in section \ref{NERML}, would return in the case of $\epsilon=\epsilon_\tau$, where $\tau$ would be the termination step. It immediately follows that at time $\tau$ we have at our disposal the {\sl best found so far} feasible solution $x^\tau$ to $(P)$ satisfying
\begin{equation}\label{Gap}
\Opt-f_*(x^\tau)\leq \Gap_\tau:=\min_{1\leq \nu\leq\tau} \epsilon_\nu
\end{equation}
(indeed, set $x^1=x_1$ and set $x^\tau=x^{\tau-1}$ when $\Gap_\tau=\Gap_{\tau-1}$ and $x^\tau=x_\tau$ otherwise). The bottom line is that instead of terminating NERML
when a given in advance accuracy $\epsilon$ is attained, we can run the algorithm for as long as we want, generating {\sl in an online fashion} the {\sl optimality gaps} $\Gap_\tau$ and feasible solutions $x^\tau$ to $(P)$ satisfying (\ref{Gap}), $\tau=1,2,...$ For experimental purposes this execution mode is much more convenient than the original one (in a single run, we get the complete ``time-accuracy'' curve instead of just one point on this curve), and this is the mode used in the experiments we are about to report.
\par
Recall that NERML is specified by a proximal setup, two control parameters $\gamma,\theta\in(0,1)$ and ``memory depth'' $m$ which should be a positive integer. In our experiments, we used the Euclidean setup (i.e., equipped the embedding space $\bR^N$ of $Y$ with the standard Euclidean norm and the distance-generating function $\omega(y)={1\over 2}y^Ty$) and $\theta=\gamma={1\over 2}$. \par
We are about to report the results of three series of experiments (all implemented in MATLAB).
\paragraph{A.}
In the first series of experiments we consider ``small'' problems ($p=512$, $r=2$, $N\in\{64,128,256,512\}$), which allows us to consider a ``wide'' range $m\in\{1,3,5,9,17,33,65,129\}$ of memory depth $m$ \footnote{In our straightforward software implementation of the algorithm, handling memory of depth $m$ requires storing in RAM up to $m+1$ $p\times p$ matrices, which makes the implementation too space-consuming when $p$ and $m$ are large; this is why in our experiments the larger $p$, the smaller is the allowed values of $m$. Note that with a more sophisticated software implementation, handling memory $m$ would require storing just $m+1$ {\sl of rank 1} $p\times p$ matrices, reducing dramatically the required space.}.  The results are presented in table \ref{table0}, where $T$ is ``physical'' running time in sec and $\Prg_m(T)$ is the {\sl progress in accuracy} in time $T$ for NERML with memory $m$, defined as the ratio $\Gap_1/\Gap_{t(T)}$,  $t(T)$ being the number of steps performed in $T$ sec.
     \par
     The structure of the data in the table is as follows. Given $p=512$, $r=2$ and a value of $N$, we generated the corresponding problem instance and then ran on this instance 1024 steps of NERML, the memory depth being set to 129, 65,...,1. As a result of these 8 runs, we got running times $T_{129}$, $T_{65}$, ..., $T_1$,  which are the values of $T$ presented in the table, and overall progresses in accuracies displayed in the column ``$\Prg_m(T)$.'' Then we ran NERML with the minimal memory 1 until the running time reached the value $T_{129}$, and recorded the progress in accuracy observed at times $T_1,...,T_{129}$, displayed in the column ``$\Prg_1(T)$.'' For example, the data displayed in the table for the smallest instance ($m=64$) say that 1024 steps of NERML with memory 129 took $\approx390$ sec, while the same number of steps with memory 1 took just $\approx 55$ sec, a $7$ times smaller time. However, the former ``time-consuming'' algorithm
     reduced the optimality gap by factor of about 8.2.e5, while the latter -- by factor of just 153, thus exhibiting about 5000 times worse progress in accuracy. We see also that even when running NERML with memory 1 for the same 390 sec as taken by the 1024-step NERML with memory 129, the progress in accuracy was ``only'' 1187 -- still by factor about 690 worse than the progress in accuracy achieved {\sl in the same time} by  NERML with memory 129. Thus, ``long memory'' can indeed be highly beneficial. This being said, we see from the table that the benefits of ``long memory'' reduce as the design dimension $N$ of the problem to which we apply NERML grows, and in order to get a ``reasonable benefit'', the memory indeed should be ``long;'' e.g., in all experiments reported in the table, NERML with memory like $9$ is only marginally better than  NERML with memory 1. The data in the table, same as the results to be reported below, taken along with the numerical experience with MDL (see the concluding comment in section \ref{sect:SCG}) allow for a ``qualified guess'' stating that {\sl in order to be highly beneficial, the memory in NERML should be of order of the design dimension of the problem at hand}.\\
     {\bf Remark:} While the numerical results reported so far seem to justify, at a qualitative level, potential benefits of re-utilizing past information, quantification of these benefits heavily depends on ``fine  structure'' and  sizes of problems in question. For example, the structure of our instances make the first order oracle for (\ref{dual_num}) pretty cheap  -- on a close inspection, a call to the oracle requires finding the leading singular vectors of a {\sl highly sparse} $p\times p$ matrix. As a result, the computational effort per step is by far dominated by the effort of solving auxiliary problems arising in NERML, and thus influence of $m$ on the duration of a step is much stronger than it could be with a more expensive first order oracle.

            \begin{table}
\caption{\label{table0}
NERML on  problem (\ref{dual_num}) with $p=512$ and $r=2$.}\vskip4pt
\centering{\small
\begin{tabular}{||c|c|c|c|c||c|c|c|c|c||}\hline\hline
$N$&$m$&$\Prg_m(T)$&$\Prg_1(T)$&$T$, sec&$N$&$m$&$\Prg_m(T)$&$\Prg_1(T)$&$T$, sec\\
\hline\hline
       &129&8.281e5&1187.2&390.1&
       &129&2.912e4&620.0&358.7\\
\cline{2-5}\cline{7-10}
        &65& 4.439e5&778.1&276.8&
        &65& 2472.5&620.0&367.0\\
\cline{2-5}\cline{7-10}
$64$ &33&2483.1&553.2&195.1&$128$   &33&217.2&564.1&207.7\\
\cline{2-5}\cline{7-10}
  &17&627.3&255.9&120.0&&17&173.2&318.4&130.5\\
  \cline{2-5}\cline{7-10}
        &9&261.1&204.0&89.1&&9&159.8&116.8&95.6\\
\cline{2-5}\cline{7-10}
        &5&184.3&154.0&71.5&&5&132.8&98.2&77.9\\
\cline{2-5}\cline{7-10}
        &3&186.2&154.0&64.7&&3&95.5&98.2&68.4\\
\cline{2-5}\cline{7-10}
        &1&152.9&154.0&55.2&&1&88.6&98.2&58.5\\
\hline\hline
            &129&5.524e4&2800.4&642.7&&129&3253.2&2850.6&944.9\\
\cline{2-5}\cline{7-10}
            &65& 1657.5&972.0&395.3&&65&350.0&1275.2&496.5\\
\cline{2-5}\cline{7-10}
$256$     &33&205.3&389.8&246.1&$512$&33&164.9&370.8&277.6\\
\cline{2-5}\cline{7-10}
       &17&159.0&267.8&147.1&&17&148.0&260.2&177.0\\
\cline{2-5}\cline{7-10}
            &9&145.9&168.7&107.5&&9&125.7&198.3&138.6\\
\cline{2-5}\cline{7-10}
            &5&132.6&125.0&87.4&&5&109.7&141.4&114.8\\
\cline{2-5}\cline{7-10}
            &3&118.2&122.0&78.0&&3&106.7&121.7&98.1\\
\cline{2-5}\cline{7-10}
            &1&87.0&97.9&66.5&&1&114.8&102.2&87.7\\
\hline\hline
\end{tabular}\\\vskip4pt
Platform: laptop PC with 2$\times$2.67GHz Intel Core i7 CPU and 8 GB RAM, Windows 7-64 OS. }

\end{table}
\paragraph{B.}
Next we apply NERML to  ``medium-size'' problems, restricting the range of $m$ to $\{1,17,33\}$  and keeping the design dimension of problems (\ref{dual_num}) at the level $N=2048$. The reported in table \ref{table1} CPU time corresponds to $t=1024$ steps of NERML. The data in the table exhibit the same phenomena as those observed on small problems.
\paragraph{C.} Finally, table  \ref{table2} displays the results obtained with 1024-step NERML with $m=1$ on ``large'' problems ($p$ up to 8196, $N$ up to 16392).
\par
We believe that the numerical results we have presented are rather encouraging. Indeed, even in the worst, in terms of progress in accuracy, of our experiments (last problem in table \ref{table1}) the optimality gap was reduced in 1024 iterations (2666 sec) by two orders of magnitude (from 0.166 to 0.002). To put this in proper perspective, note that on the same platform as the one underlying tables \ref{table1}, \ref{table2}, a {\sl single} full SVD of a 8192$\times$8192 matrix takes $>450$ sec, meaning that a proximal point algorithm applied directly to the problem of interest (\ref{pprim_num}) associated with the last line in table \ref{table1} would be able to carry out in the same 2666 sec just 6 iterations. Similarly, $\approx 3600$ sec used by NERML to solve the largest instance we have considered (last problem in table \ref{table2}, with $p=8192$ and $N=16192$, progress in accuracy by factor $\approx 1200$) allow for just 8 full SVD's of $8192\times 8192$ matrices. And of course 6 or 8 iterations  of a proximal type algorithm as applied to (\ref{pprim_num}) typically are {\sl by far} not enough to get comparable progress in accuracy.
\begin{table}
\caption{\label{table1}
NERML on medium-size problems (\ref{dual_num}), $N=2048$ and $t=1024$.}\vskip4pt
\centering
\begin{tabular}{||c|c|c|c|c|c|c||}\hline\hline
&$m$&$\Gap_1$&$\Gap_1/\Gap_{32}$
&$\Gap_1/\Gap_{128}$
&$\Gap_1/\Gap_{1024}$
&T, sec\\ \hline\hline
&33&&59.7&82.3&245.4&  390.9\\ \cline{2-2}\cline{4-7}
$p=512$, $r=4$&9&2.26e-1&23.7&76.3&227.9&  171.0\\ \cline{2-2}\cline{4-7}
&1&&14.2&60.6&168.7&  112.9\\
\hline
\hline
&33&&29.3&153.6&355.2&  757.3\\ \cline{2-2}\cline{4-7}
$p=1024$, $r=2$& 9&2.79e-1&18.6&145.3&302.9&  292.6\\\cline{2-2}\cline{4-7}
&1&&10.7&152.3&250.2&  145.9\\
\hline
\hline
&33&&16.7&27.8&206.8& 2318.8\\ \cline{2-2}\cline{4-7}
$p=2048$, $r=2$& 9&1.97e-1&16.7&30.1&132.1&  774.8\\ \cline{2-2}\cline{4-7}
&1&&16.8&33.8&147.1&  313.1\\
\hline
\hline
&33&&9.2&29.1&338.4& 7529.1\\\cline{2-2}\cline{4-7}
$p=4096$, $r=2$&9&2.55e-1&7.6&20.9&165.3& 2440.9\\ \cline{2-2}\cline{4-7}
&1&&7.8&20.6&158.8&  769.3\\
\hline
\hline
$p=8192$, $r=2$&1&1.66e-1&4.6&14.9&83.8& 2666.6\\
\hline\hline
\end{tabular}\\\vskip4pt
{\small Platform: desktop PC with   $4\times 3.40$ GHz {\em Intel Core2} CPU and 16 GB RAM, Windows 7-64 OS.}
\end{table}

\begin{table}
\caption{\label{table2}
NERML on ``large-scale'' problems (\ref{dual_num}), $m=1$}\vskip4pt
\centering
\begin{tabular}{||l|c|c|c|c|c|c||}\hline\hline
&$\Gap_1$&$\Gap_1/\Gap_{32}$&$\Gap_1/\Gap_{128}$&$\Gap_1/\Gap_{1024}$& CPU, sec\\
 \hline\hline
$p=2048$,
 $r=4$,
  $N=8192$&1.81e-1&171.2&213.8&451.4&  521.3\\ \hline
$p=4096$, $r=4$, $N=16384$
&3.74e-1&335.4&1060.8&1287.3& 1524.8\\ \hline
$p=8192$, $r=2$, $N=16384$
&2.54e-1&37.8&875.8&1183.6& 3644.0\\ \hline\hline
\end{tabular}\\\vskip4pt
{\small Platform: desktop PC with   $4\times 3.40$ GHz {\em Intel Core2} CPU and 16 GB RAM, Windows 7-64 OS.}
\end{table}

\appendix
\section{Appendix: Proofs}
We need the following  technical result originating from \cite{Teboulle93} (for proof, see \cite{Teboulle93}, or section 2.1 and Lemma A.2 of \cite{NN2013}).
  \begin{lemma} \label{AppLemma}. Let $Y$ be a nonempty closed and bounded subset of a Euclidean space $E_y$, and let $\|\cdot\|$, $\omega(\cdot)$ be the corresponding proximal setup. Let, further, $U$ be a closed convex subset of $Y$ intersecting the relative interior of $Y$, and let $p\in E_y$.\\
  {\rm (i)}  The optimization
   problem
   $$
   \min_{y\in U} h_p(y):=\left[\langle p,y\rangle + \omega(y)\right]
   $$
   has a unique solution $y_*$. This solution is fully characterized by the inclusion $y_*\in U\cap Y^o$, $Y^o=\{y\in Y:\partial \omega(y)\neq\emptyset\}$, coupled with the relation
\begin{equation}\label{addedeq1}
\langle p+\omega'(y_*),u-y_*\rangle \geq0\;\;\forall u\in U.
\end{equation}
{\rm (ii)}
When $U$ is cut off $Y$ by a system of linear inequalities $e_i(y)\leq0$, $i=1,...,m$, there exist Lagrange multipliers $\lambda_i\geq 0$ such that $\lambda_ie_i(y_*)=0,\;1\leq i\leq m$, and
\begin{equation}\label{addedeq2}
\forall u\in Y:\;\;\langle p+\omega'(y_*)+\sum_{i=1}^m\lambda_ie_i^\prime,u-y_*\rangle \geq 0.
\end{equation}
{\rm (iii)} In the situation of {\rm (ii)}, assuming $p=\xi-\omega'(y)$ for some $y\in Y^o$, we have
\begin{equation}\label{addedeq3}
\forall u\in Y: \langle \xi,y_*-u\rangle -\sum_i\lambda_ie_i(u)\leq V_y(u)-V_{y_*}(u)-V_{y}(y_*).
\end{equation}
\end{lemma}

\subsection{Proof of Proposition \ref{PropLevel}}\label{Appendix.PropLevel}
\paragraph{1$^0$.} Observe that when $t$ is a non-terminal step of the algorithm, the level set $U_t$ is a closed convex subset of $Y$ which intersects the relative interior of $Y$; indeed, by (\ref{bothcases}) and due to $\epsilon_t>0$ for a non-terminal $t$,  there exists $y\in Y$ such that $h_\tau(y)\geq \epsilon_t>\ell_t=\gamma\epsilon_t$, $\tau\in S_t$, that is, $U_t$ is cut off $Y$ by a system of constraints satisfying the Slater condition.  Denoting $S_t=\{y_{\tau_1},y_{\tau_2},...,y_{\tau_k}\}$ and invoking item (iii) of Lemma \ref{AppLemma} with $y=\widehat{y}_t$, $\xi=0$ and $e_i(\cdot)=\gamma\epsilon_t-h_{\tau_i}(\cdot)$ (so that $U_t=\{u\in Y: e_i(u)\leq 0, 1\leq i\leq k\}$, we get $y_*=y_{t+1}$ and
$$
 -\sum_i\lambda_ie_i(u)   \leq V_{\widehat{y}_t}(u)-V_{y_{t+1}}(u)-V_{\widehat{y}_t}(y_{t+1}).
$$
with some $\lambda_i\geq0$. When  $u\in U_t$, we have $e_i(u)\leq0$, that is,
\begin{equation}\label{thatis17}
\forall u\in U_t: V_{y_{t+1}}(u)\leq V_{\widehat{y}_t}(u) - V_{\widehat{y}_t}(y_{t+1}).
\end{equation}
\paragraph{2$^0$.} When $t$ starts a phase, we have $\widehat{y}_t=y_\omega=y_1$, and clearly $1\in I_t$, whence $h_\tau(\widehat{y}_t)\leq0$ for some $\tau\in I_t$ (specifically, for $\tau=1$). When $t$ does not start a phase, we have $\widehat{y}_t=y_t$ and $t\in I_t$, so that here again $h_\tau(\widehat{y}_t)\leq0$ for some $\tau\in I_t$. On the other hand, $h_\tau(y_{t+1})\geq \gamma\epsilon_t$ for all $\tau\in I_t$ due to $y_{t+1}\in U_t$. Thus, when passing from $\widehat{y}_t$ to $y_{t+1}$, at least one of $h_\tau(\cdot)$ grows by at least $\gamma\epsilon_t$. Taking into account that $h_\tau(z)=\langle g(y_\tau),y_\tau-z\rangle$ is Lipschitz continuous with constant $L[g]$ w.r.t. $\|\cdot\|$ (by (\ref{normalization})), we conclude that $\|\widehat{y}_t-y_{t+1}\|\geq\gamma\epsilon_t/L[g]$. With this in mind, (\ref{thatis17}) combines with (\ref{wehave22}) to imply that
\begin{equation}\label{thatis18}
\forall u\in U_t: V_{y_{t+1}}(u)\leq V_{\widehat{y}_t}(u) - {1\over 2}V_{\widehat{y}_t}(y_{t+1})\leq V_{\widehat{y}_t}(u)-{\gamma^2\epsilon_t^2\over2 L^2[g]}.
\end{equation}
\paragraph{3$^0$.} Let the algorithm perform phase $s$, let $t_s$ be the first step of this phase, and $r$ be another step of the phase. We claim that all level sets $U_t$, $t_s\leq t\leq r$, have a point in common, specifically, (any) $u\in \Argmax_{y\in Y}\min_{\tau\in I_r} h_\tau(y)$. Indeed, since $r$ belongs to phase $s$, we have
$$
\gamma\Delta_s<\epsilon_r=\max_{y\in Y}\min_{\tau\in I_r} h_\tau(y)=\min_{\tau\in I_r} h_\tau(u)
$$
and $\Delta_s=\epsilon_{t_s}=\max_{y\in Y}\min_{\tau\in I_{t_s}} h_\tau(y)$ (see (\ref{bothcases}) and the definition of $\Delta_s$). Besides this, $r$ belongs to phase $s$, and within a phase, sets $I_t$ extend as $t$ grows, so that $I_{t_s}\subset I_t\subset I_r$ when $t_s\leq t\leq r$, implying that
$\epsilon_{t_s}\geq\epsilon_{t_s+1}\geq...\geq\epsilon_r$. Thus, for $t\in\{t_s,t_s+1,...,r\}$ we have
$$
\min_{\tau\in I_t}h_\tau(u)\geq \min_{\tau\in I_r}h_\tau(u)\geq \gamma\Delta_s=\gamma\epsilon_{t_s}\geq
\gamma\epsilon_t,
$$
implying that $u\in U_t$.
\par
 With the just defined $u$, let us look at the quantities $v_t:=V_{\widehat{y}_t}(u)$, $t_s\leq t\leq r$. We have $v_{t_s}\leq {1\over 2}\Omega^2$ due to $\widehat{y}_{t_s}=y_\omega$ and (\ref{diameter}), and
 \[
 0\leq v_{t+1}\leq v_t -{\gamma^2\epsilon_t^2\over 2L^2[g]}\leq v_t-{\gamma^4\Delta_s^2\over2L^2[g]}\]
  when $t_s\leq t<r$ (due to (\ref{thatis18}) combined with  $\widehat{y}_t=y_t$ when $t_s<t\leq r$).  We conclude that $(r-t_s)\gamma^4\Delta_s^2\leq \Omega^2L^2[g]$. Thus, the number $T_s$ of steps of phase $s$ admits the bound
 \begin{equation}\label{Ts}
 T_s\leq {\Omega^2L^2[g]\over\gamma^4\Delta_s^2}+1\leq {2\Omega^2L^2[g]\over\gamma^4\Delta_s^2},
 \end{equation}
 where the concluding inequality follows from $\Delta_s\leq\Delta_1=\max_{y\in Y} \langle g(y_\omega),y_\omega-y\rangle \leq L[g]\Omega$, see (\ref{diameter}), combined with $\gamma\in(0,1)$.
 \paragraph{4$^0$.} Assume that MDL does not terminate in course of first $T\geq1$ steps, and let $s_*$ be the index of the phase to which the step $T$ belongs.  Then $\Delta_{s_*}>\epsilon$ (otherwise we would terminate not later than at the first step of phase $s_*$); and besides this, by construction, $\Delta_{s+1}\leq\gamma\Delta_s$ whenever phase $s+1$ takes place. Therefore
 $$
 T\leq\sum_{s=1}^{s_*} T_s\leq {2\Omega^2L^2[g]\over\gamma^4}\sum_{r=0}^{s_*-1}\Delta_{s_*-r}^{-2}
 \leq  {2\Omega^2L^2[g]\over\gamma^4\Delta_{s_*}^2}\sum_{r=0}^{s_*-1}\gamma^{2r}\leq
 {2\Omega^2L^2[g]\over\gamma^4(1-\gamma^2)\Delta_{s_*}^2}\leq {2\Omega^2L^2[g]\over\gamma^4(1-\gamma^2)\epsilon^2}.\eqno{\hbox{\qed}}
 $$

\subsection{Proof of Proposition \ref{PropNERML}}\label{Appendix.PropNERML}
 Observe that the algorithm can terminate only in the case A  of B.2.2.3, and in this case the output is indeed as claimed in Proposition. Thus, all we need to prove is the upper bound (\ref{numberofsteps}) on the number of steps before termination.

\par  1$^0$. Let us bound from above the number of steps at an arbitrary phase $s$. Assume that phase $s$ did not terminate in course of the first $T$ steps, so that $u_1,...,u_T$ are well defined. We claim that then
\begin{equation}\label{bruce22}
\|u_\tau-u_{\tau+1}\|\geq \ell_s/L[g],\,1\leq
\tau<T.
\end{equation}
Indeed, by construction $h^s_{\tau-1,m+1}(y):=\langle g(u_\tau),u_\tau-y\rangle$ is $\geq \ell_s=\gamma f_s$ when $y=u_{\tau+1}$ (due to $u_{\tau+1}\in Y_\tau$). Since $\|g(u)\|_*\leq L[g]$ for all $u\in Y$, (\ref{bruce22}) follows.
\par Now let us look at what happens with the quantities $\omega(u_\tau)$ as $\tau$ grows. By strong convexity of $\omega$ we have
$$
\omega(u_{\tau+1})-\omega(u_\tau)\geq\langle\omega^\prime(u_\tau),u_{\tau+1}-u_\tau\rangle+{1\over 2}\|u_\tau-u_{\tau+1}\|^2
$$
The first term in the right hand side is $\geq0$, since $u_\tau$ is the minimizer of $\omega(\cdot)$ over $Y^s_{\tau-1}$, while $u_{\tau+1}\in Y_\tau\subset Y^s_{\tau-1}$. The second term in the right hand side is $\geq {\ell_s^2\over 2L^2[g]}$ by (\ref{bruce22}). Thus, $\omega(u_{\tau+1})-\omega(u_\tau)\geq {\ell_s^2\over 2L^2[g]}$, whence
$
\omega(u_T)-\omega(u_1)\geq (T-1){\ell_s^2\over 2L^2[g]}=(T-1){\gamma^2f_s^2\over 2L^2[g]}.
$
Recalling the definition of $\Omega$ and that $u_1=y_\omega$, the left hand side in this inequality is $\leq{1\over 2}\Omega^2$.  It follows that whenever phase $s$ does not terminate in course of the first $T$ steps, one has
$
T\leq {\Omega^2L^2[g]\over\gamma^2f_s^2}+1,
$
that is, the total number of steps at phase $s$, provided this phase exists, is at most
$
T_s={\Omega^2L^2[g]\over\gamma^2f_s^2}+2.
$
Now, we have
$$
f_s\leq f_1=\max\limits_{y\in Y}\langle g(y_\omega),y-y_\omega\rangle \leq L[g]\max\limits_{x\in Y}\|y-y_\omega\|\leq\Omega L[g]
$$
(recall that $\|g(y)\|_*\leq L[g]$ and see (\ref{diameter})). Thus 
\[
T_s={\Omega^2L^2[g]\over\gamma^2f_s^2}+2\leq {(1+2\gamma^2)\over\gamma^2}{\Omega^2 L^2[g]\over f_s^2}
\]
for all $s$ such that $s$-th phase exists. By construction, we have $f_s\geq\epsilon$ and
$f_s\leq(\gamma +(1-\gamma)\theta) f_{s-1}$, whence the method eventually terminates (since $\gamma +(1-\gamma)\theta<1$). Assuming that the termination happens at phase
${s_*}$, we have $f_s\geq (\gamma+(1-\gamma)\theta)^{s-{s_*}}f_{{s_*}}$ when $1\leq s\leq{s_*}$, so that the total number of steps is bounded by
\bse
\lefteqn{\sum_{s=1}^{{s_*}}{(1+2\gamma^2)\over\gamma^2}{\Omega^2L^2[g]\over f_s^2}
\leq \sum_{s=1}^{{s_*}} {(1+2\gamma^2)\over\gamma^2}{\Omega^2L^2[g](\gamma+(1-\gamma)\theta)^{2({s_*}-s)}\over f_{{s_*}}^2}}\\
&\leq& \sum_{s=1}^{{s_*}} {(1+2\gamma^2)\over\gamma^2}{\Omega^2L^2[g](\gamma+(1-\gamma)\theta)^{2({s_*}-s)}\over \epsilon^2}
\leq {(1+2\gamma^2)\over\gamma^2[1-(\gamma+(1-\gamma)\theta)^2]}{\Omega^2L^2[g]\over\epsilon^2},
\ese
as claimed. \qed

\begin{thebibliography}{10}

\bibitem{Amit:etal:2007}
Y.~Amit, M.~Fink, N.~Srebro, and S.~Ullman.
\newblock Uncovering shared structures in multiclass classification.
\newblock In {\em Proceedings of the 24th international conference on Machine
  learning}, pages 17--24. ACM, 2007.

\bibitem{NERML}
A.~Ben-Tal and A.~Nemirovski.
\newblock Non-euclidean restricted memory level method for large-scale convex
  optimization.
\newblock {\em Mathematical Programming}, 102(3):407--456, 2005.

\bibitem{Bregman67}
L.~M. Bregman.
\newblock The relaxation method of finding the common point of convex sets and
  its application to the solution of problems in convex programming.
\newblock {\em USSR computational mathematics and mathematical physics},
  7(3):200--217, 1967.

\bibitem{Teboulle93}
G.~Chen and M.~Teboulle.
\newblock Convergence analysis of a proximal-like minimization algorithm using
  bregman functions.
\newblock {\em SIAM Journal on Optimization}, 3(3):538--543, 1993.

\bibitem{combettes2010}
P.~L. Combettes, D.~D{\~u}ng, and B.~C. V{\~u}.
\newblock Dualization of signal recovery problems.
\newblock {\em Set-Valued and Variational Analysis}, 18(3-4):373--404, 2010.

\bibitem{Cramer01}
K.~Crammer and Y.~Singer.
\newblock On the algorithmic implementation of multiclass kernel-based vector
  machines.
\newblock {\em The Journal of Machine Learning Research}, 2:265--292, 2002.

\bibitem{Dem:Rub:1970}
V.~Demyanov and A.~Rubinov.
\newblock {\em Approximate methods in optimization problems}, volume~32.
\newblock Elsevier Publishing Company, 1970.

\bibitem{Dunn78}
J.~C. Dunn and S.~Harshbarger.
\newblock Conditional gradient algorithms with open loop step size rules.
\newblock {\em Journal of Mathematical Analysis and Applications},
  62(2):432--444, 1978.

\bibitem{Fan11}
J.~Fan, Y.~Liao, and M.~Mincheva.
\newblock High dimensional covariance matrix estimation in approximate factor
  models.
\newblock {\em Annals of statistics}, 39(6):3320--3356, 2011.

\bibitem{Frank1956Algorithm}
M.~Frank and P.~Wolfe.
\newblock An algorithm for quadratic programming.
\newblock {\em Naval research logistics quarterly}, 3(1-2):95--110, 1956.

\bibitem{MLO}
A.~Juditsky and A.~Nemirovski.
\newblock First order methods for nonsmooth large-scale convex minimization, i:
  General purpose methods; ii: Utilizing problem's structure.
\newblock In S.~Sra, S.~Nowozin, and S.~J. Wright, editors, {\em Optimization
  for Machine Learning}, pages 121--254. Mit Press, 2011.

\bibitem{LNN}
C.~Lemar{\'e}chal, A.~Nemirovskii, and Y.~Nesterov.
\newblock New variants of bundle methods.
\newblock {\em Mathematical programming}, 69(1-3):111--147, 1995.

\bibitem{Mor62}
J.-J. Moreau.
\newblock Fonctions convexes duales et points proximaux dans un espace
  hilbertien.
\newblock {\em CR Acad. Sci. Paris S{\'e}r. A Math}, 255:2897--2899, 1962.

\bibitem{Mor65}
J.-J. Moreau.
\newblock Proximit{\'e} et dualit{\'e} dans un espace hilbertien.
\newblock {\em Bulletin de la Soci{\'e}t{\'e} math{\'e}matique de France},
  93:273--299, 1965.

\bibitem{nemirovski2004prox}
A.~Nemirovski.
\newblock Prox-method with rate of convergence $\underline{O}(1/t)$ for
  variational inequalities with lipschitz continuous monotone operators and
  smooth convex-concave saddle point problems.
\newblock {\em SIAM Journal on Optimization}, 15(1):229--251, 2004.

\bibitem{NOR}
A.~Nemirovski, S.~Onn, and U.~G. Rothblum.
\newblock Accuracy certificates for computational problems with convex
  structure.
\newblock {\em Mathematics of Operations Research}, 35(1):52--78, 2010.

\bibitem{NYu}
A.~Nemirovskii and D.~Yudin.
\newblock {\em Problem complexity and method efficiency in optimization}.
\newblock Wiley (Chichester and New York), 1983.

\bibitem{nesterov1983}
Y.~Nesterov.
\newblock A method for unconstrained convex minimization problem with the rate
  of convergence $\underline{O}(1/k^2)$.
\newblock {\em Soviet Math. Dokl.}, 27(2):372--376, 1983.

\bibitem{nesterov_book}
Y.~Nesterov.
\newblock {\em Introductory lectures on convex optimization: A basic course},
  volume~87.
\newblock Springer, 2004.

\bibitem{NesSmooth}
Y.~Nesterov.
\newblock Smooth minimization of non-smooth functions.
\newblock {\em Mathematical Programming}, 103(1):127--152, 2005.

\bibitem{Nesterov-dual}
Y.~Nesterov.
\newblock Dual extrapolation and its applications to solving variational
  inequalities and related problems.
\newblock {\em Mathematical Programming}, 109(2-3):319--344, 2007.

\bibitem{Nesterov-primal}
Y.~Nesterov.
\newblock Primal-dual subgradient methods for convex problems.
\newblock {\em Mathematical programming}, 120(1):221--259, 2009.

\bibitem{NN2013}
Y.~Nesterov and A.~Nemirovski.
\newblock Some first order algorithms for $\ell_1$/nuclear norm minimization.
\newblock {\em Acta Numerica}, pages 1--67, 2013.

\bibitem{Pshe:1994}
B.~N. Pshenichnyi and Y.~M. Danilin.
\newblock {\em Numerical methods in extremal problems}.
\newblock Mir Publishers (Moscow), 1978.

\bibitem{Yos64}
K.~Yosida.
\newblock {\em Functional analysis}.
\newblock Springer Verlag, Berlin, 1964.

\end{thebibliography}
\end{document}